\documentclass[12pt]{article}
\usepackage{amsmath}

\usepackage{amsfonts}
\usepackage{graphicx}
\usepackage{latexsym}
\usepackage{amscd,amsfonts}
\usepackage{amssymb,amsmath}

\newtheorem{theorem}{\bf Theorem}
\newtheorem{lemma}{\bf Lemma}
\newtheorem{example}{\bf Example}
\newtheorem{remark}{\bf Remark}
\newtheorem{corollary}{\bf Corollary}

\newtheorem{definition}{{\bf Definition}}

\begin{document}

\title{
\bf\large Dimension results for inhomogeneous Moran set constructions}
\author{Mark Holland, Yiwei Zhang}
\date{\today}
\maketitle



\begin{abstract}
We compute the Hausdorff, upper box and packing dimensions for
certain inhomogeneous Moran set constructions. These constructions
are beyond the classical theory of iterated function systems, as
different nonlinear contraction transformations are applied at each
step. Moreover, we also allow the contractions to be weakly
conformal and consider situations where the contraction rates have
an infimum of zero. In addition, the basic sets of the construction
are allowed to have a complicated topology such as having fractal
boundaries. Using techniques from thermodynamic formalism we
calculate the fractal dimension of the limit set of the
construction. As a main application we consider dimension results
for stochastic inhomogeneous Moran set constructions, where chaotic
dynamical systems are used to control the contraction factors at
each step of the construction.
\end{abstract}

\noindent{\bf Keywords:} Fractal dimension, inhomogeneous Moran
sets.

\noindent{\bf Subject Classification:} 28A78, 28A80, 37A05, 37F35.

\section{Introduction}
A systematic study on the classical theory of iterated function
systems (IFS) has been developed in the pioneering work of Moran
\cite{Moran} and Bowen \cite{Bowen}, and has been successfully
applied in the study of dimension theory (e.g. Bowen's formalism for
$C^{1+\epsilon}$ repellers \cite{ChenPesin,Falconer97,Przytycki09}).
However, most scenarios require the iterated function system IFS
to be \emph{conformal},
and step independent.
In this paper we go beyond these classical settings, and
consider inhomogeneous Moran set constructions. The main
difficulties encountered on estimating the fractal dimension for these constructions
are as follows. Firstly, the nonlinear contractions in the IFS are
step dependent. Secondly, these contractions are allowed to be
weakly conformal (in ways that we will make precise). We also
allow the basic sets of the construction to have wild topological
properties (such as fractal boundaries), and permit arbitrary
placement of the basic sets, subject to these sets being separated.
To study inhomogeneous Moran set constructions, we combine various
approaches such as those considered in \cite{FWW, Hua00, Ma01} and
\cite{Barreira, ChenPesin, Pe97,PesinWeiss}. Our aim is to form a
unified approach in the computation of fractal dimension for such
inhomogeneous constructions.

To obtain concrete results on the fractal dimensions such as
Hausdorff, upper-box and packing dimension we introduce the main
geometrical hypotheses in Section \ref{sec:definitions}. These
include assumptions on the degree of nonlinearity permitted on the
contractions, and control on the placement of the basic sets in
terms of their separation (rather than their precise location).
Within this section we also introduce the mechanism of symbolic
codings used to describe the basic sets of the construction. In
particular, when the IFS is affine or one dimensional
cookie-cutter-like, our dimension results on inhomogeneous Moran
sets coincide with the results obtained in
\cite{FWW,Hua00,Ma01,Wen01}. In our weakly conformal case, we permit
no specific control on the distortion or smoothness of the
contraction maps except for continuity. Instead we concentrate on
the cardinality of the Moran covering as well as the existence of a
\emph{Gibbs-like} measure. We also consider constructions defined
on sub-symbolic spaces. In particular, we consider
sub-spaces formed by placing restrictions on the sequence of
admissible words, for example by introducing a transition matrix. We
study the corresponding fractal dimension when the sequence of
admissible words is restricted, see Section \ref{subsec:subspace}.
These constructions can be viewed as generalized versions of
\emph{graph directed Markov systems} (see \cite{MU03}).

The main dimension results are presented in Section
\ref{sec:mainresults}, where we determine the fractal dimension of a
limit set $F$ in terms of a sequence of pre-dimensions $s_k$. The
pre-dimension sequence depends on the first $k$ steps of the
construction, and for non-linear constructions we take $s_k$ to be
the zero of a corresponding pressure equation $P_k(s\Phi_k)=0$, with
a defined potential $\Phi_k$, see Section \ref{subsec:predim}. For
nonlinear constructions of inhomogeneous Moran sets, our approach
extends the theory developed in \cite{PesinWeiss}, where they
primarily control the geometry using a single vector (of contraction
constants) with a finite number of components. In our case we work
with a countable sequence of vectors, and the geometry of the
construction is controlled using this vector sequence, see Section
\ref{subsec:vectors}. Moreover we consider scenarios where the
infimum of the contraction vector components is equal to zero, and
comment on situations where the supremum of the contraction vector
components equals 1. For example, we believe our techniques will
extend to inhomogeneous constructions generated by nonlinear cookie
cutters with parabolic fixed points. In the context of IFS having
parabolic fixed points, see \cite{GelfertRams, Urbanski1,
Urbanski2}.

As another novelty, we also consider stochastic constructions of
inhomogeneous Moran sets and give corresponding dimension results.
This is discussed in Section \ref{sec:applications}. For such
constructions, we use a stationary stochastic process to generate
the $k$-step contraction rates, for example by taking a time series
of observations on an ergodic transformation (see \cite{Walters82}).
This approach appears to be new, at least relative to classical
stochastic constructions mentioned in \cite{Falconer, Wen01}. This
gives an alternative approach for constructing random fractals using
ergodic and statistical properties of dynamical systems. We study
the typical (almost sure) fractal dimension, and further
investigations might include studying the largest/smallest
dimensions that can arise (e.g. utilizing ideas from ergodic
optimization theory \cite{Jenkinson}). We further consider
stochastic constructions where the infimum of the contraction vector
components is equal to zero, and where the corresponding
supremum equals 1.

The formal proofs of the dimension results are presented in Section
\ref{sec:proof}, with background on dimension theory and
thermodynamic formalism presented in Section \ref{sec:appendix}.

\section{Geometric and symbolic constructions}\label{sec:definitions}
\subsection{Symbolic spaces for inhomogeneous Moran set constructions}\label{subsec:definitons}
We define the following symbolic space. For a sequence of positive
integers $\{n_{k}\}_{k\geq 1}$ and any $k\in \mathbb{N}$, let
\begin{equation}\label{equ_symbol}
D_{k}=\{(i_{1},i_{2},\cdots,i_{k}); ~~~1\leq i_{j}\leq n_{j}, 1\leq
j\leq k\}\quad\textrm{with}\quad D_{0}=\emptyset,
\end{equation}
and define
\begin{equation}\label{equ_symbol2}
D=\bigcup_{k=0}^{\infty}D_{k}.
\end{equation}
The set $D_k$ contains all words of length $k$. The collection
$D$ is a countable collection of level sets.
\begin{definition}
Given a map $f:\mathbb{R}^{d}\rightarrow \mathbb{R}^{d}$, we define
the class $\Im$ such that if $f\in\Im$ then
\begin{enumerate}
  \item[(H1)] There exists a compact forward invariant set $A$, such that $f(A)\subset A$;
  \item[(H2)] For any compact set $B\subset A$, $\mathrm{diam}(f^{n}(B))\rightarrow0$ as $n\rightarrow\infty.$
\end{enumerate}
\end{definition}
We remark that for any $f\in\Im$, then
$\bigcap_{n=1}^{\infty}f^{n}(A)$ is a singleton. If $f_{1},f_{2}\in
\Im$ share the same forward invariant set $A$, then both
$f_{2}\circ f_{1}$ and $f_{1}\circ f_{2}\in \Im$. We say that a map
$f$ is \emph{contracting} if there exists a
$0<c<1$ such that for all $x,y\in \mathbb{R}^{d}$, $d(f(x),f(y))\leq
c\cdot d(x,y)$. If $f$ is contracting then $f\in\Im$, but the converse
need not be true.

\begin{definition}
A family of compact sets is called {\bf basic sets}
$\Omega=\{\triangle_{\omega}\subset\mathbb{R}^{d},\,\omega\in D\}$,
if this family of sets satisfies: $\lim_{k\to\infty}\max_{\omega\in
D_{k}}\operatorname{diam}(\triangle_{\omega})=0.$
\end{definition}

Based on $D$ and the class $\Im$ of maps, we consider the
following \emph{Moran structure conditions} (MSC) for a
class of sets $\Omega=\{\triangle_{\omega},\,\omega\in D\}$, where
$\triangle_{\omega}\subset\mathbb{R}^d$ and $\omega=(i_1,i_2,\ldots,
i_k)$ is a finite \emph{word} in $D$. Given words $\omega,\omega'\in
D$ we define $\omega*\omega'$ as the concatenation of the two words (when this is
still defined in $D$).

\begin{definition}\label{def_mor_like}
Given a basic set $\triangle\subset\mathbb{R}^d$ and a sequence of
contractions $\{f_{j,i}\Im:i\leq n_j, j\geq 1\}$ we say that
$\Omega=\{\triangle_{\omega},\,\omega\in D\}$ satisfies (MSC) with
respect to $D$ if the following hold.
\begin{enumerate}
\item[(A1)] Suppose $k\geq 1$, $\omega\in D_{k-1}$ and $\omega*j\in D_{k}$ (for $1\leq j\leq n_k$).
Then elements of $\triangle_{\omega*j}$ are completely determined by
elements of $\triangle_{\omega}$ and the vector of maps
$\Xi_{k}=(f_{k,1},f_{k,2},\ldots,f_{k,n_k})\in\Im$, i.e.,
$\triangle_{\omega*j}\subseteq \triangle_{\omega}$ and
$\triangle_{\omega*j}=f_{k,j}(\triangle_{\omega})$.
Moreover if $\omega=(i_1,i_2,\ldots, i_k)$, then we define
$f_{\omega}=f_{1, i_1}\circ f_{2,i_{2}}\circ\cdots\circ f_{k,i_{k}}$, so that
$f_{\omega}(\Delta)=\Delta_{\omega}$.
\item[(A2)] The strong separation condition holds: given any $k$ and $\omega,\omega'\in D_k$ with $\omega\neq\omega'$ then
$$\triangle_{\omega}\cap\triangle_{\omega'}=\emptyset.$$
\end{enumerate}
\end{definition}
For a given $\Omega$, we define
\begin{equation}\label{moranF}
E_k=\bigcup_{\omega\in D_k}\triangle_{\omega},\;
F=\bigcap_{k\in\mathbb{N}}E_k.
\end{equation}
The set $F$ is a compact set, and by the strong separation condition
(A2) is totally disconnected. So far we have made no assumptions on
the topology of the basic set $\triangle$, nor on the sets
$\triangle_{\omega}$ ($\omega\in D$) other than these sets being
compact. In particular they need not to be connected, and their
boundaries could be fractal. It is sufficient for our purposes to
work with a weaker version of (A2), and we say that the \emph{weak
separation condition holds} if
\begin{enumerate}
\item[(A2')] For any $\omega,\omega'\in D$ with $\omega\neq\omega'$:
$$\{\triangle_{\omega}\cap\triangle_{\omega'}\}\cap F=\emptyset.$$
\end{enumerate}

\begin{figure}[h]
  \centering
  \includegraphics[width=8cm]{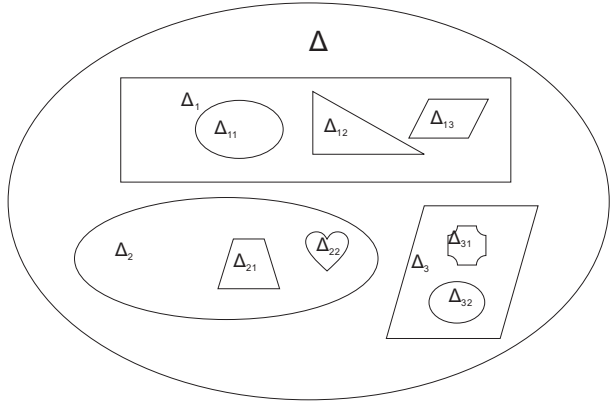}\\
  \caption[Simple geometric construction]{Schematic representation of the geometric construction of a general
Moran set $F$.}
  \label{fig_fractal}
\end{figure}

\begin{definition}\label{def_mor_like2}
Given $F$ as in equation \eqref{moranF}, we call $F$ a {\bf
generalized Moran set} (GMS) if $F$ satisfies (A1) and (A2').
\end{definition}
See Fig \ref{fig_fractal} for the geometrical interpretation.

Now define the set
\[
D^{*}=\{(i_{1},i_{2},\cdots, i_{k},\cdots): 1\leq i_{j}\leq n_{j},\,
j\geq1\}.
\]
This set consists of infinite strings, and any $\omega\in D^*$ has the representation
$\omega=(i_1,i_2,\ldots).$ Given $\omega\in D^*$, we write $C_{(i_1,\ldots i_k)}(\omega)\subset D^*$
as the $k$-length cylinder set. Given $D^{*}$ and $F$, there is a canonical projection map
$\mathcal{X}:D^{*}\to F$ which assigns to each
$\omega=(i_n)_{n=1}^{\infty}$ the point $x\in F$ given by
$\bigcap_{k}\triangle_{(i_1,\ldots i_k)}$.

We can turn $D^{*}$ into a metric space by assigning the distance
function $d(\omega,\omega')$ to points $\omega',~\omega\in D^{*}$ as
follows:
$$n(\omega,\omega'):=
          \min\{i|\omega_{j}=\omega_{j}'~~\mbox{for}~0<j<i~\mbox{but}~\omega_{i}\neq
\omega_{i}'\},\mbox{~if~} \omega\neq \omega',$$ and
$n(\omega,\omega):=\infty$. For given $p_i<n^{-1}_{i}$ we set
$d(\omega,\omega')=\prod_{i=1}^{n(\omega,\omega')}p_{i}.$ Then
$(D^{*},d)$ is a compact metric space. It is easy to see that
$D^{*}$ is a generalization of traditional symbolic space, since if
$n_{k}=p$ is a constant, then $D^{*}=\Sigma_{p}^{+}$, where
$\Sigma_{p}^{+}=\{1,\ldots,p\}^{\mathbb{N}}$. When $n_k$ is not
constant the shift map $\sigma$ on $D^*$ does not preserve $D^*$ in
general. We consider a sequence of symbolic spaces that can be
thought as approximations to $D^*$. These symbol spaces are
generated from the sets $D_k$.

\begin{definition}\label{def_shiftmap}
Given $D=\bigcup_k D_k$, the symbol space $[D_k]$ is defined as the
set of infinite strings, with indices corresponding to elements of
$D_k$. That is
$$[D_k]=\{\mathbf{w}=(\omega_i)_{i=1}^{\infty}=(\omega_1,\omega_2\,\ldots,),\,\omega_i\in D_k\}.$$
The associated shift map $\sigma_{k}:[D_k]\to [D_k]$ is defined by:
\begin{equation*}
  \sigma_{k}(\omega_1,\omega_2\,\ldots,)=(\omega_2\,\omega_3,\ldots,),\quad\textrm{with}\;(\omega_i)_{i=1}^{\infty}\in [D_k].
\end{equation*}
\end{definition}
Notice that $[D_k]$ is isomorphic to $\Sigma^{+}_{p_k}$ with
$p_k=\mathrm{card}(D_k)$.

\subsection{Sub-spaces of symbolic constructions}\label{subsec:subspace}
So far we have considered all admissible collections of words in $D=\cup_k D_k$. Instead, we can consider
subsets of words $Q_k\subset D_k$, with $Q=\cup_k Q_k\subset D$. If $\{A^{(k)}\}$ is a sequence of (transition) matrices, having
entries in $\{0,1\}$ then admissible words in $Q$ may be characterized in terms of products of these matrices.
In particular we can write
\begin{equation}
Q_k=\{\omega=(i_1,\ldots, i_k)\in D_k:\,
A^{(1)}_{i_1i_2}A^{(2)}_{i_2i_3}\cdots
A^{(k-1)}_{i_{k-1}i_k}=1\},\quad Q=\bigcup_{k} Q_k.
\end{equation}
Thus with $Q$ and $Q_k$ in place of $D$, resp. $D_k$, we can produce constructions in analogy to those considered
in Definitions \ref{def_mor_like} and \ref{def_mor_like2}, but now for the class of sets
$\Omega(Q)=\{\Delta_{\omega}\subset\mathbb{R}^d,\omega\in Q\}$. The corresponding limit set $F$ defined by
\begin{equation}\label{moran-subF}
E_k=\bigcup_{\omega\in Q_k}\triangle_{\omega},\quad
F=\bigcap_{k\in\mathbb{N}}E_k,
\end{equation}
will be referred to as a generalized Moran set associated to $Q$. We define
\[
Q^{*}=\{(i_{1},i_{2},\cdots, i_{k},\cdots):\,A^{(k)}_{i_{k}i_{k+1}}=1, k\geq 1\},
\]
and given $\omega\in Q^*$, we write $C_{(i_1,\ldots i_k)}(\omega)\subset Q^*$
as the $k$-length cylinder set. There is again a canonical projection map
$\mathcal{X}:D^{*}\to F$ that takes $\omega=(i_n)_{n=1}^{\infty}$ to the $x\in F$ given by
$\bigcap_{k}\triangle_{(i_1,\ldots i_k)}$. We again can turn $Q^*$ into a metric space
(using the metric inherited from that of $D^*$), and we define the symbol space $[Q_k]$ in direct
analogy to $[D_k]$.

Since $Q$ can be quite general, we will mainly consider the case where the transition matrices
$A^{(k)}:=A$ are fixed $p\times p$ matrices (and hence $n_k=p$ for each $k$). We can then find the
fractal dimension of $F$ in terms of the (spectral) properties of $A$, and in terms of the contraction
vector sequence $\Xi_k$ as defined in condition (A1).

\subsection{Conformal constructions and constructions bounded via upper/lower estimating vectors}\label{subsec:vectors}
To obtain explicit estimates on the Hausdorff dimension of $F$, some
restrictions on the basic sets $\triangle_\omega$ are required. In
particular we require control on the diameter of
$\triangle_{\omega}$ with respect to the level set $D_k$ that
$\omega$ belongs to. In particular we require that their diameters
shrink exponentially fast with $k$. We also require control the
geometry of $\triangle_{\omega}$ via a technical condition
restricting the number of $\triangle_{\omega}$ (of a certain
size-scale) that can intersect with a given ball
$B(x,r)\in\mathbb{R}^d$ where $x\in F$. For self similar
constructions, control on the geometry is specified in
\cite{PesinWeiss} by use of lower, and upper estimating vectors. We
adapt these methods for the non-self similar constructions. Let
$\overline\Psi=\{\overline\Psi^{(k)},\,k\in\mathbb{N}\}$ denote a countable
collection of vectors $\overline\Psi^{(k)}$ with
$$\overline\Psi^{(k)}=(\Psi^{(k)}(\omega))_{\omega\in D_k}.$$
Here $\omega$ has the representation as some $(i_1,\ldots, i_k)\in
D_k$. 
Given $\omega\in D_k$, we assume that there is a sequence of
constants $c_{i_1},\ldots c_{i_k}$ such that
$$\Psi^{(k)}(\omega)=c_{i_1}^{(1)}c_{i_2}^{(2)}\ldots c_{i_k}^{(k)}=\prod_{j=1}^{k}
c^{(j)}_{i_j}.$$
For notational simplicity we sometimes write
$\Psi^{(k)}_{\omega}:=\Psi^{(k)}(\omega)$. In relation to the sequence $\overline\Psi^{(k)}$ we define $\tilde{\Xi}_k$ to be the $k$-step vector sequence:
$$\tilde{\Xi}_k=(c^{(k)}_1,c^{(k)}_2,\ldots,c^{(k)}_{n_k}).$$
For example, if the $k$-step vector $\Xi_k=(f_{k,i})_{i=1}^{n_k}$ consists of affine maps, each with contraction rate
$c^{(k)}_{i_k}$, then a natural choice for $\tilde{\Xi}_k$ would be the vector sequence of corresponding contraction ratios.

\begin{definition}[Basic vectors] \label{defbasic}
The collection of vectors
$\overline\Psi=\{\overline\Psi^{(k)}\}_{k=1}^{\infty}$ is called a
{\bf basic collection} of vectors if for all $k\geq 1$ and all
$\omega=(i_1,\ldots, i_k)\in D_k$, the sequence $\Psi^{(k)}(\omega)$
satisfies
\begin{equation}\label{equ_basic}
\sup_{k\in\mathbb{N},1\leq j\leq n_{k}}c^{(k)}_j<1.
\end{equation}
\end{definition}

\begin{definition}[UE vectors]\label{defupper}
A basic collection of vectors
$\overline\Psi=\{\overline\Psi^{(k)}\}_{k=1}^{\infty}$ is called an
{\bf upper estimating} (UE) collection of vectors if for any $k$ and
$\omega\in D_{k}$:
\[
\operatorname{diam}(\triangle_{\omega})\leq C\Psi^{(k)}_{\omega},
\]
and the constant $C>0$ is independent of $\omega$ and $k$.
\end{definition}

To get bounds on the Hausdorff dimension we require further control
of the geometry of each $\Delta_{\omega}$. We introduce two
definitions: the first is that of conformality, while the second
introduces the notion of lower-estimating vectors for a geometric
construction.

\begin{definition}[Conformal vectors]\label{defconformal} Given a basic collection of vectors
$\overline\Psi=\{\overline\Psi^{(k)}\}_{k=1}^{\infty}$, we say that
a symbolic construction $\{\triangle_{\omega}\}$ is {\bf conformal}
(w.r.t. $\overline\Psi$) if $\exists\, C>0$, such that for each $k\geq
1$, $\omega\in D_k$, $\exists x\in\Delta_{\omega}$:
\begin{equation}
B\left(x,\frac{1}{C}\Psi^{(k)}_{\omega}\right)\subset\triangle_{\omega}\subset
B\left(x,C\Psi^{(k)}_{\omega}\right).
\end{equation}
\end{definition}

The following geometric constraint is formulated in terms of Moran coverings
which we describe as follows, see also
\cite{PesinWeiss}. Given a set $F$, and for any $x\in F$, choose
the $\omega\in D^{*}$ for which $\mathcal{X}(\omega)=x$. By the
separation condition, $\omega$ is unique. Suppose $0<r<1$ is fixed and let
$\overline\Psi$ be a basic sequence of vectors. Let $n(x)$ be the
unique positive integer of the such that
\[
\Psi^{(n(x))}_{\omega}>r~\mbox{and}~\Psi^{(n(x)+1)}_{\omega}\leq r.
\]
If $C(\omega)$ is the corresponding $n(x)$-length cylinder set, we write
$\triangle(x):=\mathcal{X}(C(\omega))$. For $x,y\in F$, either
$\triangle(x)=\triangle(y)$ or $\triangle(x)\cap\triangle(y)=\emptyset$.
The corresponding (disjoint) collection of sets we denote by
$\{\triangle^{(j)}\}$, where $F\subset\cup_j\triangle^{(j)}$, and this forms
the \emph{Moran covering} of the set $F$.

Consider the open ball $B(x,r)$ of the radius $r$ centered at the
point $x\in F$, and let $N(x,r)$ denote the cardinality of the subset
of $\{\triangle^{(j)}\}$ that have non-empty intersection with $B(x,r)$. We
have the following definition.

\begin{definition}[LE vectors]\label{def_regular}
If there exists a constant $M$ such that the above $N(x,r)<M$ for
all $x\in F$, then we say the collection of vectors $\overline\Psi$
is {\bf lower estimating} (LE).
\end{definition}

In the special case where the vector $\overline\Psi$ has the property that $\Psi^{k}_{\omega}=\Psi^{(k)}_{\omega'}$
for all $\omega,\omega'\in D_k$, then we call the construction \emph{homogeneous} if such a vector
is both (UE) and (LE). The corresponding limit set $F$ is called homogeneous, otherwise
in all other cases the construction (and limit set) is  \emph{inhomogeneous}.

\subsection*{Pre-dimension sequences}\label{subsec:predim}

For MSCs arising from non-linear constructions, we determine the dimension of the Moran
set $F$ from a sequence of pre-dimensions $s_k$.
These $s_k$ will be prescribed to be the zeros
of a functional equation involving the topological pressure. We make
this precise as follows. Consider a sequence of pressure functions
$P_k$ (for $k\in\mathbb{N}$), and a sequence of potentials $\Phi_k$ defined as
follows. Suppose that $\overline\Psi$ is prescribed, and consider the symbolic space $[D_k]$ together with the shift
map $\sigma_k:[D_k]\to [D_k].$ For
$\mathbf{w}=(\omega_1,\omega_2,\ldots)\in [D_k]$ let
$\Phi_{k,s}(\mathbf{w}):=s\log\Psi^{(k)}(\omega_1)$. This function
can be extended to a function on $F_k$ via
$\Phi_{k,s}(x)=s\log\Psi^{(k)}(\omega_1)$, where
$\mathcal{X}(\mathbf{w})=x$. We define the corresponding pressure
function $P_k:\mathrm{Lip}(F_k)\to\mathbb{R}$ by
\begin{equation}\label{k-pressure}
P_k(\Phi_{k,s})=\lim_{n\rightarrow\infty}\frac{1}{n}\log
\left(\sum_{(\omega_1,\ldots,\omega_n)}
\underset{x\in\Delta_{(\omega_1,\ldots,\omega_n)}}{\inf}
\exp\left\{S_{n}(\Phi_{k,s}(x))\right\}\right),
\end{equation}
where $\Delta_{\omega_1,\ldots,\omega_n}=f_{\omega_1,\ldots,\omega_n}(\Delta)$,
and $\omega_i\in D_k$. Now we consider the sequence
$s_k$,  where $s_k$ is the value of $s$ which solves
$P_k(\Phi_{k,s})=0$. In particular we consider the (lim)-inf and
(lim)-sup of this sequence. We define:
\begin{equation}\label{s-limit}
s^{*}:=\lim\sup s_k,\quad\textrm{and}\quad s_{*}:=\lim\inf s_k.
\end{equation}
The main focus of this paper is to consider when $s^{*}$ is the
upper-box dimension of $F$, and when $s_{*}$ is the Hausdorff
dimension of $F$.

To obtain dimension estimates for $F$ in terms of zeros of the pressure
function we need to assume the existence of a \emph{Gibbs-like} measure
on $F$ as follows:
\begin{enumerate}
\item[(A3)] Given $\beta>0$, there exists a measure $m_{\Psi}$ supported on $F$, and $L>0$ such that
for all $k\geq 1, \omega\in D_k$,
\begin{equation}
\frac{L^{-1}}{\sum_{\omega'\in D_k}(\Psi^{(k)}(\omega'))^{\beta}}
\leq\frac{m_{\Psi}(\Delta_{\omega})}{(\Psi^{(k)}(\omega))^{\beta}}\leq
\frac{L}{\sum_{\omega'\in D_k}(\Psi^{(k)}(\omega'))^{\beta}}
\end{equation}
\end{enumerate}
For a range of applications hypothesis (A3) can be verified. For IFS
defined by expanding maps, then (A3) typically follows from bounded
distortion estimates, see Section \ref{sec:applications}. Without
(A3), assumptions (A1) and (A2) will not ensure that
$\dim_H(F)=s_*$.

\section{Statement of main results}\label{sec:mainresults}

For general Moran set constructions
we now compute (or estimate bounds) on the Hausdorff, upper-box and packing
dimensions based on the existence of a countable sequence $\overline\Psi$ of upper and lower estimating vectors.
We will assume geometrical assumptions (A1), (A2) and existence of a Gibbs-like measure (A3). Applications fitting
these geometrical models will be discussed in Section \ref{sec:applications}. The constant $c_*$ will also be of
importance, where we define
\begin{equation}\label{equ_inf}
c_{*}:=\inf_{k\in \mathbb{N},1\leq j\leq n_{k}}c^{(k)}_j.
\end{equation}
We will distinguish between cases where $c_*>0$ and $c_*=0$.
\begin{theorem}\label{main.theorem}
Consider a MSC with $F$ a GMS. Suppose
$\overline\Psi=\{\overline\Psi^{(k)},\,k\in\mathbb{N}\}$ is a basic sequence of
vectors which satisfy the (UE), (LE) properties, and suppose that
there exists a Gibbs-like measure $m_{\Psi}$ satisfying (A3). Assume
further that $c_{*}>0$, where $c_*$ is defined in equation
\eqref{equ_inf}. Then
\begin{enumerate}
\item $\dim_{H}F=\mathrm{dim}_{H}m_{\Psi}= s_{*}$.
\item $\dim_{P}F,\,\overline{\dim_{B}}F\leq s^{*}.$
\end{enumerate}
If instead $F$ satisfies the conformality condition, as in
Definition \ref{defconformal} then
$$\dim_{P}F=\overline{\dim_{B}}F= s^{*}.$$
\end{theorem}
We remark that under assumption \eqref{equ_inf}, the existence of a
conformal vector
implies the (LE) property; see the proof of Lemma \ref{Conformality
implies l-estimating}. However, the converse does not hold.
Under the assumption of a vector being lower estimating, and the
construction non-conformal then we only obtain the inequality
$\overline{\dim_{B}}F\leq s^{*}.$ So far, we do not have an explicit
construction of a set $F$ for which the inequality is strict.

Suppose now that $c_*=0$.
Then we have to impose conditions on how fast the $\Psi^{(k)}_{\omega}$
decay to get corresponding results as stated in Theorem \ref{main.theorem}.
For fixed $k$, we denote
\begin{equation}\label{equation-m-d}
M_{k}:=\max\limits_{\omega\in D_{k}}\Psi^{(k)}_{\omega},
\;d_{k}:=\min\limits_{1\leq j\leq
n_{k+1}}c^{(k+1)}_{j}.
\end{equation}
We have the following result
\begin{theorem}\label{main.theorem.2}
Consider a MSC with $F$ a GMS.  Suppose
$\overline\Psi=\{\overline\Psi^{(k)},\,k\in\mathbb{N}\}$ is a sequence of vectors
which satisfy the (UE), (LE) properties, and suppose that there
exists a Gibbs-like measure $m_{\Psi}$ satisfying (A3). Furthermore
assume that $c_*=0$, and
\begin{equation}\label{equ_limcondition1}
  \lim_{k\rightarrow\infty}\frac{\log d_{k}}{\log M_{k}}=0,
\end{equation}
then $\dim_{H}F=\dim_H(m_{\Psi})=s_{*}$ and
$\dim_{P}F,\,\overline{\dim_{B}}F\leq s^{*}$, where $d_{k}, M_{k}$
are defined in \eqref{equation-m-d}.
If instead $F$ satisfies the conformality condition, as in
Definition \ref{defconformal} then
$$\dim_{P}F=\overline{\dim_{B}}F= s^{*}.$$
\end{theorem}
It is possible to impose alternative conditions on the vectors
$\Psi^{(k)}_{\omega}$ where \eqref{equ_inf} holds. We consider
the following conditions, suppose
\begin{gather}
M=\sup_{k\geq 1}n_k<\infty\label{inf-zero2a}\\
0<\inf_{k}\max_{1\leq j\leq n_k}c^{(k)}_{j}
\leq\sup_{k}\max_{1\leq j\leq n_k}c^{(k)}_{j}<1.\label{inf-zero2b}
\end{gather}
For example, equation \eqref{equ_limcondition1} can be satisfied for
a homogeneous construction having $c^{(k)}_j=c_k$, for $1\leq j\leq n_k$,
and $\inf_k c_j=0$. However for this example, equation \eqref{inf-zero2b}
will fail. An example that satisfies \eqref{inf-zero2b}, but not \eqref{equ_limcondition1} would be a construction
with vector $\tilde\Xi_k=(1/4, 1/4, (1/4)^k).$ The following theorem holds.
\begin{theorem}\label{main.theorem3}
Consider a MSC with $F$ a GMS.  Suppose
$\overline\Psi=\{\overline\Psi^{(k)},\,k\in\mathbb{N}\}$ is a \emph{basic} sequence
of vectors which satisfy the (UE), (LE) properties, and suppose that
there exists a Gibbs-like measure $m_{\Psi}$ satisfying (A3).
Moreover, suppose that equations \eqref{inf-zero2a},
\eqref{inf-zero2b} hold with $c_*=0$. Then
$\dim_{H}F=\dim_H(m_{\Psi})=s_{*}$. If instead, $\overline\Psi$ is a
\emph{conformal vector} then $\dim_{P}F=\overline{\dim_{B}}F=
s^{*}.$
\end{theorem}
When $\sup_{k,j} c^{(k)}_j=1$ and/or when $\sup_k n_k=\infty$ then it is possible to give constructions where
$\dim_H(F)\neq\lim\inf s_k$, and/or $\overline{\dim_{B}}(F)\neq\lim\sup s_k$, see \cite{Hua00}.
For Moran set constructions
modelled by subsets of symbolic spaces then corresponding results hold. We state the following
corollary (whose proof follows step by step from the proofs of Theorems  \ref{main.theorem}, \ref{main.theorem.2}
and \ref{main.theorem3}).

\begin{corollary}\label{coro-subshift}
Suppose that $F$ is a GMS generated by a sub-symbolic space
$Q_k\subset D_k$, with $n_k=p$ fixed, and allowed words modelled by
a (fixed) transition matrix $A$. Relative to the space $Q$, suppose
$\overline\Psi=\{\overline\Psi^{(k)}\}$ is a sequence of vectors
which satisfy the (UE), (LE) properties. Furthermore suppose that
relative to the space $Q$ there exists a Gibbs-like  measure
$m_{\Psi}$ satisfying (A3). Then the conclusions of Theorems
\ref{main.theorem}, \ref{main.theorem.2} and \ref{main.theorem3}
remain valid.
\end{corollary}

\section{Applications}\label{sec:applications}
We consider applications of Theorems \ref{main.theorem},
\ref{main.theorem.2} and \ref{main.theorem3} to a range of examples.
We first consider step dependent IFS, and then explore Moran set constructions with
stochastic vectors.

\subsection{Iterated function systems}\label{sec_reference}
In this section we consider IFS defined by sequences of expanding
maps.

Suppose that we are given a basic set
$\triangle\subset\mathbb{R}^{d}$ and
$\Omega=\{\triangle_{w}\in\mathbb{R}^{d}:\omega\in D\}$ satisfies the conditions of
MSC as stated in Definition \ref{def_mor_like}.
Based on these geometrical constructions, we consider a family of maps $\{T_{i,j}\}$ defined
in the following way.

$T_{j,i_{j}}:\triangle_{i_1,\ldots,i_j}\to\triangle_{i_1,\ldots, i_{j-1}},\forall
i_j=1,\cdots,n_{k}$ satisfies the following assumptions:
\begin{description}
  \item[(IFS1):] For $j\geq 1$ and $1\leq i_j\leq n_j$,
$T_{j,i_{j}}:\triangle_{i_1,\ldots,i_j,}\to\triangle_{i_1,\ldots, i_{j-1}}$ is a full-branch $C^{1+\alpha}$ diffeomorphism.
In particular $T_{j,i_{j}}(\triangle_{i_1,\ldots,i_j,})=\triangle_{i_1,\ldots, i_{j-1}},\,\forall
i_j=1,\cdots,n_{k}$, and
the derivative $DT_{j,i_{j}}$ is $\alpha-$ H\"{o}lder continuous, i.e., there exists a constant $C:=C_{j,i_{j}}$
such that $||D T_{j,i_{j}}(x)-DT_{j,i_{j}}(y)||\leq C||x-y||^{\alpha}$.
 \item[(IFS2):] There exists a $\beta:=\beta_{j,i_{j}}>1$ such
  that $||T_{j,i_{j}}(x)-T_{j,i_{j}}(y)||\geq\beta||x-y||,~\forall
  x,y\in\triangle_{i_1,\ldots, i_{j}}$;
\end{description}

We take $\Xi_k=(f_{k,1},\ldots, f_{k,n_k})$ to be the
vector of contractions associated to the inverse branches of
$(T_{k,1},\ldots,T_{k,n_{k}})$ at the $k-$th step. For
$\omega=(i_1,\ldots, i_k)\in D_k$ we have
$\triangle_{\omega}=f_{\omega}(\triangle)$ where
$f_{\omega}=f_{1,i_1}\circ\cdots,\circ f_{k,i_k}$. We state the following
result

\begin{corollary}\label{coro-ifs}
For a family of expanding diffeomorphisms $\{T_i\}$, let
$\{\Xi_{k}\}_{k=1}^{\infty}$ be the vector sequence of contractions
associated to the inverse branches. Consider a GMS,
$F$ associated to this $\{\Xi_k\}$. Assume that the $C^{1+\alpha}$
(distance)-expansivity of $\{\Xi_{k}\}_{k=1}^{\infty}$ is uniformly
bounded, i.e., the sequence $\{\beta_{j,i_{j}}\}$ is uniformly
bounded away from 1, the sequence $\{\det(Df_{j,i_{j}})\}$ is
uniformly bounded away from zero, and the sequence of H\"{o}lder
constants $\{C_{j,i_{j}}\}$ is uniformly bounded. Then
$$\dim_{H}F= s_{*},\quad\dim_{P}F=\overline{\dim_{B}}F= s^{*},$$
where $s_{*}$ and $s^{*}$ are defined in equation \eqref{s-limit}.
\end{corollary}

Before giving the proof consider the example where
$f_{i}$ are similarity contractions
and the basic sets $\triangle_{\omega}$ as intervals (or balls) in
$\mathbb{R}^d$, see \cite{Hua00}. We show how the corresponding dimension
estimates are obtained by assuming (A1), (A2), and checking (A3).
The problem can be reduced to taking a sequence of vectors $\tilde\Xi_k$ (associated to $\Xi_k$) given by
\begin{equation}
\tilde\Xi_{k}=(c^{(k)}_1,c^{(k)}_{2},\ldots,c^{(k)}_{n_k}),
\end{equation}
where the $c^{(k)}_{i}=Df_{k,i}\mid\Delta$ are positive constants.
Assuming (A1) and (A2) there is a similarity transformation $f_{\omega}$ taking $\Delta$ to $\Delta_{\omega}$. Moreover,
suppose $k\geq 1$, $\omega\in D_{k-1}$ and $\omega*j\in D_{k}$
(for $1\leq j\leq n_k$).
Then $\triangle_{\omega*j}\subset \triangle_{\omega}$, and
$$\frac{|\Delta_{\omega*j}|}{|\Delta_{\omega}|}=c^{(k)}_{j}.$$
The corresponding pre-dimension sequences $\{s_k\}$ satisfy the equations
\begin{equation}\label{eq:linear-construction}
\prod_{i=1}^{k}\sum_{j=1}^{n_j}(c^{(i)}_{j})^{s_k}=\sum_{\omega\in D_k}(\mathrm{diam}f_{\omega}(\Delta))^{s_k}=1.
\end{equation}
These equations are equivalent to solving $P(s_{k}\Phi_{k}(x))=0$,
where $\Phi_k(x)=s_k\log\Psi^{(k)}(\mathbf{w}^{(1)})$, $\forall
x\in\Delta_{\omega_1,\ldots,\omega_n}$, and $P(\cdot)$ is defined in equation \eqref{k-pressure}.
The corresponding Gibb-like measure $m_{\Psi}$ can be made taken as
the weak limit of the sequence of measures $m_k$, where each $m_k$ is defined on
$\omega\in D_{\ell}, \ell\leq k$ as follows:
\begin{equation*}
m_{k}(\Delta_{\omega})=\sum_{i_{\ell+1},\ldots,i_k}\frac{(c_{1,i_1}c_{2,i_2}\ldots c_{k,i_k})^{\beta}}
{\prod_{j=1}^{k}\sum_{i=1}^{n_k}c_{j,i}^{\beta}},\quad\omega=(i_1,\ldots i_{\ell}).
\end{equation*}
By linearity of the construction we have $m_{k}(\Delta_{\omega})=m_{\ell}(
\Delta_{\omega})$. The estimates are uniform in $k$, and hence (A3) holds
when taking a weak limit of $\{m_k\}$. We therefore obtain by Theorem
\ref{main.theorem}
\begin{equation}
\dim_H(F)=s_*,\quad \overline{\dim_B}(F)=\dim_{P}(E)=s^{*},
\end{equation}
where
\begin{equation}
s_{*}=\liminf_{k\to\infty}s_k,\quad s^{*}=\limsup_{k\to\infty} s_k.
\end{equation}

\medskip

\noindent{\em Proof of Corollary \ref{coro-ifs}:} The key calculation in the nonlinear setting is to use
bounded distortion. We show that the construction can be modelled
by a basic and conformal vector sequence $\overline{\Psi}$. Furthermore
we check that (A3) holds.

First of all, we claim that there exists
$D>0$, independent of $k$ such that for all $x,y\in\triangle$ and
$\omega\in D_k$
\begin{equation}\label{equ_bounded distortion}
\frac{1}{D}\leq
\frac{|\det(Df_{\omega}(x))|}{|\det(Df_{\omega}(y))|}\leq
D.
\end{equation}
The proof of the distortion result is based on the chain rule, for
the same iterated function system at each level; see
\cite{Falconer,Przytycki09}. More precisely, we have:

\begin{align*}
&\left|\log|\det Df_{\omega}(x)|-\log|\det
Df_{\omega}(y)|\right|\nonumber\\
=&\sum_{j=1}^{k}\left|\log|\det Df_{j,i_{j}}(f_{\omega\mid
j}(x))|-\log|\det
Df_{j,i_{j}}(f_{\omega\mid j}(y))|\right|\\
\leq & \sum_{j=1}^{k}C_{1}\left|\det D_{j,i_{j}}(f_{\omega\mid
j}(x))-\det
D_{j,i_{j}}(f_{\omega\mid j}(y))\right|\\
\leq &\sum_{j=1}^{k}C_{2}||Df_{j,i_{j}}(f_{\omega\mid
j})(x)-Df_{j,i_{j}}
(f_{\omega\mid j})(y)||\\
\leq &
\sum_{j=1}^{k}C_{3}||f_{\omega\mid j}(x)-f_{\omega\mid j}(y)||^{\alpha}\\
\leq &
C_{3}\sum_{j=1}^{k}\beta^{-j\alpha}||x-y||^{\alpha}\leq\frac{C_{3}\beta^{-\alpha}}{1-\beta^{-\alpha}}||x-y||^{\alpha},
\end{align*}
where for $j\leq k$, $\omega=(i_1,\ldots, i_k)\mid j$ corresponds to the word
$(1_1,\ldots,i_j)$.
Due to the uniform bounded distortion, these constants
$C_{i},i=1,2,3$ and $\beta$ are independent of the choice of $k$,
which implies \eqref{equ_bounded distortion}.

From this bounded distortion property \eqref{equ_bounded
distortion}, we can directly construct a collection of vectors
$\overline\Psi$ and verify the conformality and (A2). More
precisely, for any fixed $x\in\triangle$, let
$\Psi_{\omega}=\sup_{x\in\triangle_{\omega}}|\det Df_{\omega}(x)|, \forall
\omega=\omega\in D_{k}$. Then, for all $y\neq x\in\triangle$,
we have
$$\frac{1}{D}\leq\frac{|\det Df_{\omega}(y)|}{\Psi^{(k)}_{\omega}}\leq D.$$
Thus by the expanding and distortion properties of $\{T_i\}$,
the vector sequence $\overline\Psi$ is basic and conformal. To check
assumption $(A3)$ we take $m_{\Psi}$ as weak limit of
measures $m_k$, where each $m_k$ is defined on
$\omega\in D_{\ell}, \ell\leq k$ as follows:
\begin{equation*}
m_{k}(\Delta_{\omega})=\frac{\left(\mathrm{diam}(\triangle_{\omega})\right)^{\beta}}
{\sum_{\omega\in D_k}\left(\mathrm{diam}(\triangle_{\omega})\right)^{\beta}}.
\end{equation*}
This is in complete analogy to the linear construction considered for
similarity transformations. A computation using bounded distortion, see
\cite[Prop 2.7]{Ma01}, implies that $m_{\Psi}$ satisfies (A3). The corresponding
results on the dimension follow from Theorem \ref{main.theorem}. \hfill $\Box$

Corollary \ref{coro-ifs} extends the results of \cite{Ma01} to
higher dimensions, and to situations where the basic sets have fractal boundaries.
The results also apply when taking instead complex conformal holomorphic
expanding maps on the Riemann sphere $\overline{\mathbb{C}}$. In this case we
let $\Psi_{\omega}^{(k)}:=\max_{x\in\Delta_{\omega}}|\operatorname{Arg}(f'_{\omega}(x))|$).

So far we have assumed the sequence of vectors to be basic. The authors conjecture
that this assumption can be relaxed, and the results extend to the scenario where
the class of maps $\{T_i\}$ are \emph{non-uniformly expanding}. An example would
include the parabolic-fixed point family of
maps $T_i:[0,1]\to\mathbb{R}$, $\alpha_i\in(0,1)$ given by
\begin{equation}\label{intermittent}
T_{i}(x)=
\begin{cases}
x(1+3x^{\alpha_i}) &\textrm{if}\,x\in[0,1/2],\\
3(1-x) &\textrm{if}\,x\in(1/2,1].
\end{cases}
\end{equation}
For $\alpha=\alpha_i$ fixed, and potential $\phi(x)=s\log T'(x)$
the corresponding pressure function is no longer
analytic in $s$. There is a critical value $s=s_c$ for which
the pressure function undergoes a phase transition (corresponding to
derivative singularity). For all $s>s_c$, the pressure function is zero. However it can be
shown that $\dim_H(F)=s_c=\inf\{s:P(s\phi)=0\},$ see \cite{GelfertRams,
Urbanski1, Urbanski2}.
For inhomogeneous Moran set constructions generated by a sequence of maps $T_i$.
The authors conjecture that for a sequence of maps $T_i$, each
having a parabolic fixed point (with parabolic index $\alpha_i$) the corresponding
dimension is given by $\dim_H(F)=s_*$, with $s_{*}=\liminf_k s_k$, and
$s_k=\inf\{s:P_k(s\Phi_k)=0\}$.

\subsection{Stochastic Moran set constructions}\label{sec:stochastic}
In this section we consider Moran set constructions based on stochastic vector models.
Given $\omega=(i_1,\ldots, i_k)\in D_k$, we assume the constants $c^{(j)}_{i_j}(\omega)$
that constitute the vector $\overline\Psi^{(k)}$ are generated by a stationary stochastic process, such as an
ergodic transformation.

\subsubsection*{Homogeneous-stochastic Moran set constructions}
The homogeneous construction is perhaps the simplest example of a
MSC. A natural exploration is to consider ways of generating the
limit set $F$ via stochastic sequences of contractions. For example,
we consider the vector $\overline\Psi$ generated stochastically via
chaotic maps in the following sense: Let $(T,M,\mu)$ be a measure
preserving system, where $T:M\to M$ is a map preserving an ergodic
measure $\mu$. Given a test function (observable) $\phi:M\to [0,1]$
and initial condition $x\in M$, we let
$\Psi^{(k)}_{\omega}=\prod_{j=1}^{k}\phi(T^j(x))$ for any $\omega\in
D_k.$ We assume that $n_k=q$ is fixed, and the conditions of
Definition \ref{def_mor_like} apply. In this case the vector $\Xi_k$
consists of $q$ components each with value $\phi(T^k(x))$. Thus the
limit set $F$ (and hence its dimension) depends on the initial value
$x\in M$. In this section we primarily investigate the Hausdorff
dimension of $F$, and it's dependency on $x$. The results are
obtained by using methods in ergodic theory.
\begin{theorem}\label{theorem_ergodic}
Suppose that $(T,M,\mu)$ is an ergodic system, and suppose that
$\phi:M\to [0,1)$ is such that $\log\phi\in L^1(\mu)$ with $\int\log\phi<0$.
Suppose further that $F$ is the homogeneous
GMS arising from a MSC with a basic vector
$\Psi^{(k)}_{\omega}=\prod_{j=1}^{k}\phi(T^j(x))$ that is both (LE) and (UE). Assuming (A3),
then for $\mu$-a.e. $x\in M$
\[
\dim_{H}(F)=\dim_{P}(F)=\overline{\dim_{B}}(F)=\frac{-\log
q}{\int\log\phi d\mu}.
\]
\end{theorem}

\noindent{\em Proof:} When $\inf_{x\in M}\phi(x)>0$, Theorem \ref{main.theorem} implies that
\[
\dim_{H}(F)=s_{*},
\]
where
\[
s_{*}=\liminf_{k\to\infty}\left(\frac{\log\left(\prod_{j=1}^{k}\phi\circ
T^{j}(x)\right)}{-k\log q}\right)^{-1}.
\]
The Birkhoff Ergodic Theorem implies that  $\mu$-a.e. $x\in M:$
\[
\lim_{k\rightarrow\infty}\frac{1}{k}\log\left(\prod_{j=1}^{k}\phi\circ
T^{j}(x)\right)
=\lim_{k\rightarrow\infty}\frac{1}{k}\sum_{j=1}^{k}\log\phi\circ
T^{j}(x)=\int\log\phi d\mu.
\]
Now consider the case where $\inf_{x\in M}\phi(x)=0$. Since
$\inf\limits_{k\in\mathbb{N},1\leq j\leq
n_{k}}c_{j}^{(k)}=0$ we need to show that equation
\eqref{equ_limcondition1} applies for $\mu$-typical orbits, and then we apply
Theorem \ref{main.theorem.2}.

Using the notation of equation
\eqref{equ_limcondition1} we have $d_{k}=\phi(T^{k+1}(x))$ and
$M_{k}=\prod_{i=1}^{k}\phi(T^{i}(x)).$ Hence,
\begin{equation}\label{equ_propsition2}
\frac{\log d_{k}}{\log
M_{k}}=\frac{\log\phi(T^{k+1}(x))}{\sum_{i=1}^{k}\log\phi(T^{i}(x))}=\frac{k^{-1}\log\phi(T^{k+1}(x))}{k^{-1}\sum_{i=1}^{k}\log\phi(T^{i}(x))}.
\end{equation}
Again, by the ergodic theorem,
$k^{-1}\sum_{i=1}^{k}\log\phi(T^{i}(x))\to\int\log\phi\,d\mu\neq 0$.
To show $k^{-1}\log\phi(T^{k+1}(x))\to 0$ (for $\mu$-a.e. $x\in M$),
let $a_k=k$, $b_k=\log(\phi(T^k(x))$ and $S_k:=\frac{1}{a_k}\sum_{j=1}^{k}b_k$.
Then
\begin{equation*}
\frac{a_{k+1}}{a_{k}}S_{k+1}-S_k=\frac{b_{k+1}}{a_{k}},
\end{equation*}
and taking limits on both sides implies $\lim_{k\to\infty}b_{k+1}/a_k=0$.
Hence
$\log d_{k}/\log M_{k}\to 0$ for $\mu$-a.e. $x\in M$, proving the result.

\hfill $\Box$\\

\subsubsection*{Inhomogeneous-stochastic Moran set constructions}\label{sec.stoch-multi}
Consider a family of maps $\{(T_i,M,\mu_i)\}_{i=1}^{q}$ with
$T_i:M\to M$ ($M$ compact), and each $T_i$ preserves an ergodic
measure $\mu_i$ with density in $L^{p}$ for some $p>1$. Given $\mathbf{x}\in M^{q}$, we can
generate a limit set $F$ via a MSC in the following way. Take
continuous functions $\phi_i:M\to[0,1]$, and suppose that the basic
vector $\Psi^{(k)}_{\omega}$ has the form:
\begin{equation}
\Psi^{(k)}_{\omega}=\prod_{j=1}^{k}\phi_{i_j}(T^{j}_{i_j}(x_{i_j})):\;
\omega=(i_1,\ldots, i_k),\;\mathbf{x}=(x_1,\ldots,x_q).
\end{equation}
We have the following result.
\begin{theorem}\label{theor_ergodic}
Suppose that $\{T_i,M,\mu_i\}_{i=1}^{q}$ form a mixing system
(i.e., each measure $\mu_{i}$ is mixing w.r.t. $T_{i}$) and each
$\phi_i:M\to [0,1)$ is positive H\"older continuous with $\int
\log\phi_i\,d\mu_i<0$. Suppose that $F$ is a GMS arising
from a MSC with a basic vector $\overline\Psi$ generated via the
vectors $\tilde\Xi_k=(\phi_1(T^{k}_1(x_1)),\ldots,\phi_t(T^{k}_t(x_t))).$
We also assume that the basic vectors satisfy the (UE),(LE)
properties and $(A3)$ condition. Then for $\mathbf{\mu}$-a.e. $\mathbf{x}\in
M^{q}$,
\[
\dim_{H}(F)=\dim_{P}(F)=\overline{\dim_{B}}(F)=s_*,
\]
where $s_*$ is the unique solution of the functional equation:
\begin{equation}\label{eq:dimension_sol}
I_s:=\int_{M^{q}}\log\left\{\sum_{i=1}^{q}\phi_i(x_i)^{s_*}\right\}d\mathbf{\mu}=0,\;\textrm{where}\;
\mathbf{\mu}=\mu_1\times\mu_2\dots\times\mu_q.
\end{equation}
\end{theorem}

For classical stochastic (and statistically self-similar) constructions, e.g. those described in
\cite{Falconer}, they instead consider the contraction ratios $|\Delta_{\omega*j}|/|\Delta_{\omega}|:=C_j(\omega)$
as independent and identically distributed random variables. i.e. For each $j$, $\{C_{j}(\omega),\omega\in D_k\}$
are identically distributed and independent, although for fixed
$\omega$ the set of random variables $\{C_{j}(\omega),j\leq n_{k+1}\}$ need not be independent. The corresponding Hausdorff
dimension $s$ satisfies the expectation equation $E(\sum_{j=1}^{q}C_{j}^s)=1$, which is \emph{not} equivalent to
\eqref{eq:dimension_sol}. The result of \cite{Falconer} is proved using a combination of martingale and potential
theoretic methods.
Consider inhomogeneous Moran set constructions, where the contraction ratios
$C_j(\omega)$, $\omega\in D_k$ are governed by probability distributions that vary with step $k$. Then
under suitable geometric constraints, see \cite{Wen01} the Hausdorff
dimension of the corresponding limit set is given by $\dim_{H}(F)=s_*$, where
$$s_{*}=\liminf_{k\geq 1} s_k,\quad E\left(\sum_{(i_1,\ldots, i_k)}\prod_{j=1}^{k}(C_{i_j})^s\right)=1.$$
Random symbolic constructions are also included in \cite{PesinWeiss}, and these include constructions with
random vectors. They do not specifically generate the stochasticity using chaotic maps, and in their case they
obtain only the inequality $\dim_H(F)\geq s$, where $s$ satisfies the equation $\sum_{i=1}^{q}\exp\{s\int\log\phi_i\,d\mu_i\}=1.$
By a reverse Minkowski inequality this is consistent with the equality we obtain in \eqref{eq:dimension_sol}.

\medskip

\noindent{\em Proof of Theorem \ref{theor_ergodic}:} We first consider the case where
$\inf_{i}\inf_{x_{i}}\phi_{i}(x_{i})>0$. Since the set $F$ results
from a MSC, conditions (A1)-(A3) hold and it is implicit that the
$\phi_i$ are contractions.
The corresponding contraction vector is given by
$\tilde\Xi_k=(\phi_1(T^{k}_1(x_1)),\ldots,\phi_q(T^{k}_q(x_q)).$ It
suffices to compute the pre-dimensions $s_k$ and calculate the limit
$\liminf s_k$. We
have:
\begin{equation}
\sum_{\omega\in
D_k}\left(\prod_{j=1}^{k}\phi_{i_j}(T^{j}_{i_j}(x_{i_j}))\right)^{s_k}=1.
\end{equation}
A simple application of the binomial theorem implies that this
expression is equivalent to:
\begin{equation}\label{eq.mulit.pre-dimension}
\prod_{j=1}^{k}\left(\sum_{i=1}^{q}\{\phi_{i}(T^{j}_{i}(x_i)\}^{s_k}\right)=1,
\end{equation}
and so
\begin{equation}\label{eq:multi.erg}
\sum_{j=1}^{k}\log\left(\sum_{i=1}^{q}\{\phi_{i}(T^{j}_{i}(x_{i})\}^{s_k}\right)=0.
\end{equation}
Now for \emph{fixed} $s$, and by the ergodic theorem, we have
\begin{equation}\label{eq:multi.erg2}
\lim_{k\to\infty}\frac{1}{k}\sum_{j=1}^{k}\log\left(\sum_{i=1}^{q}\{\phi_{i}(T^{j}_{i}(x_{i})\}^{s}\right)=
\int_{M^{q}}\log\left\{\sum_{i=1}^q\phi_i(x_i)^{s}\right\}d\mathbf{\mu}.
\end{equation}
In the above, we have used the fact that the product system is ergodic. This is true provided
each $\mu_i$ is mixing, \cite{Walters82}.
Clearly, the value $s_*$ which is the solution of
\eqref{eq:dimension_sol} gives the right hand side of
\eqref{eq:multi.erg2} as zero. By monotonicity of $I_s$,
the value of $s_*$ is unique. We now justify that $s_*=\liminf s_k$
by showing that for large $k$, $s_k=s_{*}+o(1)$. For finite (but
large) $k$, we have
\begin{equation}
\sum_{j=1}^{k}\log\left(\sum_{i=1}^{q}\{\phi_{i}(T^{j}_{i}(x_{i})\}^{s}\right)=
k\left(\int_{M^q}\log\left\{\sum_{i=1}^q\phi_i(x_i)^{s}\right\}d\mathbf{\mu}+o(1)\right).
\end{equation}
By continuity of $\phi_i$, it follows that $\forall\,\epsilon>0$, the exists a $K$
such that $\forall k\geq K$, we can choose $s_k$ with
$|s_*-s_k|<\delta$ and $s_k$ satisfying \eqref{eq:multi.erg}. Hence
$s_*=\liminf s_k$.

Suppose now that $\inf\limits_{1\leq i\leq
q}\{\inf_{x_i}\phi_i(x_i)\}=0$, but $\int\phi_i\,d\mu_i\neq 0$. We now have $c_*=0$, see equation \eqref{equ_inf}.
Therefore, we need to show that equation
\eqref{equ_limcondition1} applies for $\mathbf{\mu}$-typical orbits.
If so, then Theorem \ref{main.theorem.2} will establish the
corresponding result. Proceeding, and using the notation of equation
\eqref{equ_limcondition1} we have that
\begin{equation}
\begin{split}
d_{k} &=\textrm{min}_{1\leq i\leq q}\{\phi_i(T^{k+1}_i(x_i))\},\\
M_{k}
&=\textrm{max}_{\omega}\left\{\prod_{j=1}^{k}\phi_{i_j}(T^{j}_{i_j}(x_{i_j}))\right\},\quad\omega=(i_1,\ldots
i_k).
\end{split}
\end{equation}
We now show the following
\begin{lemma}
Under the hypothesis of Theorem \ref{theor_ergodic}, we have for $\mathbf{\mu}$-a.e. $\mathbf{x}\in M^t$
\begin{equation*}
\lim_{k\to\infty}\frac{\log d_k}{\log M_k}=0
\end{equation*}
\end{lemma}

\noindent{\em Proof:}
We first notice that there is a constant $\lambda<1$ such that
\begin{align*}
&M_{k} \leq (\sup_{i}\sup_{x_{i}}\phi_{i}(x_i))^k\leq\lambda^k,\\
\implies &\log M_k \leq k\log\lambda<0,
\end{align*}
and
\begin{equation*}
-\log d_k=\textrm{max}_{1\leq i\leq q}\{-\log \phi_i(T^{k+1}_i(x_i))\}>0.\\
\end{equation*}
Together these imply that
\begin{equation}\label{est_bc}
\frac{\log d_k}{\log M_k}
\leq\frac{\textrm{max}_{1\leq i\leq
q}\{-\log\phi_i(T^{k+1}_i(x_i))\}}{-k\log(\sup_{i}\sup_{x_{i}}\phi_i(x_{i}))}
\leq
\frac{\textrm{max}_{1\leq i\leq
q}\{-\log\phi_i(T^{k+1}_i(x_i))\}}{-k\log\lambda} .
\end{equation}
We have to show that for $\mu_i$-a.e. $x_i$, the right hand term of
equation \eqref{est_bc} goes to zero. We use a Borel-Cantelli
argument as follows. Let
$$A^{(i)}_{k}=\{x_i\in M:\,\phi_i(T^{k}_{i}(x_i))\leq \lambda^{\sqrt{k}}\}.$$
If $\mathbf{x}\in M^q$ is such that $x_i\not\in A^{(i)}_{k}$ (for
each component $x_{i}$), then $\phi_i(T_{i}^k(x_i))<\lambda^{\sqrt{k}}$ and $\log d_k> \sqrt{k+1}\log\lambda.$
By invariance of $\mu$, and the fact that $\mu\in L^p$ we have by H\"older's inequality
$$\mu(A^{(i)}_k)\leq C(\mathrm{Leb}\{x\in M:\phi_i(x)<\lambda^{\sqrt{k}}\})^{\frac{(p-1)}{p}},$$
where $C$ depends only on $\phi_i$.
Moreover, by H\"older continuity of $\phi_i$, there is a constant $\gamma>0$ such that $\mu(A^{(i)}_k)\leq C\lambda^{\gamma\sqrt{k}}$,
and hence
$$\sum_{k=1}^{\infty}\mu_i(A^{(i)}_{k})\leq C\sum_{k=1}^{\infty} \lambda^{\gamma\sqrt{k}}<\infty.$$ Therefore by
the Borel-Cantelli Lemma that for $\mu$-a.e. $\mathbf{x}\in M^{q}$,
\begin{equation}\label{est2_bc}
\frac{\log d_k}{\log M_k}\leq\lambda^{\sqrt{k}},\quad (\textrm{eventually as
$k\to\infty$}).
\end{equation}
\hfill $\Box$

We consider two further examples which can be easily generalized to other scenarios.
\begin{example}
{\rm Take $q=2$, and suppose $\phi_1(x):=\phi(x)$ is non-constant with
$\inf_{x\in M}\phi(x)=0$. Suppose also that $\phi_2(x)=\lambda<1/2$ (constant). In this example
$M_k\geq\lambda^k$, and it need not be true that $\log d_k/\log M_k\to 0$. However in this example
Theorem \ref{main.theorem3} applies, and the corresponding dimensions are given by equation
\eqref{eq:dimension_sol}.}
\end{example}
\begin{example} {\rm Consider the case where $\sup_{x\in M}\phi_i(x)=1$ (for at least one $i\leq q$).
Take $q=2$, and for a given function $\phi(x):M\to [0,1]$, and
constant $\lambda<1$ let $\phi_1(x)=\phi(x)$,
$\phi_2(y)=\lambda(1-\phi(y))$. We will set $x=y$ and this
dependency is to ensure that $\phi_1+\phi_2< 1$ for all steps of the
construction. We also take $T_1=T_2$. If $\int\log\phi\,d\mu<0$ then
an application of the ergodic theorem tells us that for $\mu$ a.e.
$x$, the corresponding upper estimating vectors
$\Psi^{(n)}_{\omega}$ are basic. If $\int\log\phi\, d\mu=0$ then the
upper-estimating vector need not be basic, and an explicit example
would be to take $T(x)$ an interval map with a parabolic index
greater than 1, i.e. a map of the form given in equation
\eqref{intermittent}. If the upper-estimating vector is not basic
then $\overline{\dim}_B(F)$ will depend on both the placement of the
basic sets and the pre-dimension sequence $s_k$.}
\end{example}

\subsection{Stochastic Moran set constructions defined on subspaces of symbolic spaces}
We consider a stochastic Moran set construction in the setting
of Section \ref{sec.stoch-multi} with
$n_k=p$ fixed, and $A^{(k)}$ a fixed $q\times q$ matrix. Hence
\begin{equation}
Q_k=\{\omega\in D_k: A_{i_1i_2}A_{i_2i_3}\cdots A_{i_{k-1}i_k}=1\},\quad Q=\bigcup_k Q_k.
\end{equation}
We consider the (dimension) properties of the limit set $F$ defined by
$$F=\bigcap_{n\geq 1}\bigcup_{\omega\in Q_k}\Delta_{\omega}.$$
Our aim is to obtain a corresponding formula for the fractal dimension of $F$ in terms of the limiting sequence $s_k$
(defined
using the full word sequence $D_k$), and the spectral properties of $A$. In the following we let $\rho(A)$ denote
the spectral radius of a matrix $A$.
\begin{theorem}\label{theor_subshift}
Suppose that $\{T_i,M,\mu_i\}_{i=1}^{q}$ form an ergodic system
(i.e., each measure $\mu_{i}$ is ergodic w.r.t. $T_{i}$) and each
$\phi_i:M\to [0,1]$ is positive H\"older continuous with $\int
-\log\phi_i\,d\mu_i<\infty$. Suppose that the GMS, $F$ arises from a
MSC with a basic vector $\overline\Psi$ generated via the vector sequence
$\tilde\Xi_k=(\phi_1(T^{k}_1(x_1)),\ldots,\phi_q(T^{k}_q(x_q))),$ and
fixed transition matrix $A^{(k)}=A$. We also assume that the basic
vectors satisfy the (UE), (LE) properties and the $(A3)$ condition. Then
for $\mathbf{\mu}$-a.e. $\mathbf{x}\in M^{q}$
\[
\dim_{H}(F),\dim_{P}(F),\overline{\dim_{B}}(F)\leq s_*.
\]
Here $s_*$ is the unique solution of the functional equation:
\begin{equation}\label{eq:dimension_sol_shift}
\int_{M^{q}}\log\left\{\rho\left(A^{T}\Phi(\mathbf{x},s_*)\right)\right\}d\mathbf{\mu}=0,\;\textrm{where}\;\mathbf{\mu}=\mu_1\times\mu_2\cdots\times\mu_q
\end{equation}
and $\Phi(\mathbf{x},s)$ is the diagonal matrix $\mathrm{diag}(\phi_1(x_1))^s,\ldots,\phi_q(x_q)^s),$.
\end{theorem}
\begin{example}
Suppose admissible elements in $Q$ are characterized by a $p\times p$ transition matrix $A$ taking values in $\{0,1\}$, so that an element  $\omega=(i_1,i_2,\ldots)\subset Q$ is admissible if $A_{i_ji_{j+1}}=1$. Suppose we take contractions generated
by a family of similarities with constant contraction rate $\alpha_i$, $1\leq i\leq q$ (and independent of the Moran construction step $k$). Then a straightforward calculation, see \cite{PesinWeiss} implies that
\begin{equation}
\dim_H(F)=\dim_B(F)=\dim_P(F)=\tilde{s},\;\textrm{with}\;\rho(A^T\mathrm{diag}(\alpha^{\tilde{s}_1},\ldots,\alpha^{\tilde{s}_q})=1.
\end{equation}
\end{example}
\begin{remark}
Notice that we only obtain inequality in Theorem \ref{theor_subshift}. If $A$ is a constant matrix of 1s then we obtain Theorem \ref{theor_ergodic}
as before.  If the spectral radius of $A$ is equal to 1 then $s^*$ is equal to zero, and hence the corresponding dimensions are zero.
\end{remark}
\begin{remark}
The proofs would adapt easily to more general situations where $A^{(k)}$ varies with $k$. However,
explicit bounds on the fractal dimensions in terms of the spectral properties of $A^{(k)}$ are
perhaps less tractable.
\end{remark}
\begin{remark}
If instead we have a homogeneous construction with stochastic vector
$\Psi^{(k)}_{\omega}=\prod_{j=1}^{k}\phi(T^j(x))$. Then for $\mu$-a.e. $x\in M$
\[
\dim_{H}(F),\,\dim_{P}(F),\,\overline{\dim_{B}}(F)= -\frac{\log\rho(A)}{\int_M\log\phi(x)\,d\mu}.
\]
\end{remark}

\noindent{\em Proof:} Following the proof of Theorem
\ref{theor_ergodic}, the corresponding equation that replaces
equation \eqref{eq.mulit.pre-dimension} is the following:
\begin{equation}\label{eq.predim.spec}
\rho\left(\prod_{i=1}^{k}A^T\Phi_{i}(\mathbf{x},s_k)\right)=1,
\end{equation}
where $\Phi_i(\mathbf{x},s)=(\mathrm{diag}(\phi_1(T_{1}^i(x_1))^s,\ldots,\phi_p(T_{p}^i(x_p))^s).$
We can reduce this equation to the inequality $s_k\leq \tilde{s}_k$, where
\begin{equation}
\sum_{i=1}^{k}\rho\left(A^T\Phi_{i}(\mathbf{x},\tilde{s}_k)\right)=0.
\end{equation}
This follows from $\rho(XY)\leq \rho(X)\rho(Y)$ (for matrices $X,Y$), and also from the fact that the left hand side of
equation \eqref{eq.predim.spec} is monotonically decreasing in $s$. We can now take limits in $k$ as in the proof of Theorem
\ref{theor_ergodic} and hence obtain $\tilde{s}_*=\lim\inf_k s_k\leq s_*$, where $s_*$ satisfies
equation \eqref{eq:dimension_sol_shift} as stated in the Theorem. In the case where $\inf\phi_i=0$
equation \eqref{est_bc} still applies for $d_k$ and $M_k$ when restricted to admissible words in $Q$. Hence, if
the corresponding vectors satisfy properties (UE), (LE) and there exists a Gibbs measure satisfying (A3) then
$dim_{H}(F),\dim_{P}(F),\overline{\dim_{B}}(F)\leq s_*$ as required.

\section{Proof of main results}\label{sec:proof}
\subsection{Proof of Theorem \ref{main.theorem}}\label{sec:proof1}

The proof of Theorem \ref{main.theorem} is given in several steps.
We first obtain an estimate on how the measure $m_{\Psi}$ scales on
balls of radius $r$, as $r\to 0$. A second step is to show
$\dim_{H}(F)=s_{*}$ using the
existence of a basic vector sequence with the (LE), (UE) properties.

\begin{lemma}\label{lemma1}
Suppose $\overline\Psi=\{\overline\Psi^{(k)}\}$ is a \emph{basic}
sequence of vectors which satisfy the (LE) property and (A3). Take a GMS, $F$ and $x\in F$. Then for any open ball $B(x,r)$,
$(0<r<1)$ and any $\epsilon>0$, there exists a Gibbs measure
$m_{\Psi}$ such that
\begin{equation}\label{eq.meas.reg}
m_{\Psi}(B(x,r))\leq Cr^{s_*-\epsilon}.
\end{equation}
The constant $C>0$ is independent of $r$.
\end{lemma}

\noindent{\em Proof:} For $\omega\in D_k$ and for any
$\beta<s_{*}:=\lim\inf s_{k}$, by property (A3), we have
\begin{equation}
m_{\Psi}(\triangle_{\omega})\leq L_{1}(\Psi^{(k)}_{\omega})^{\beta}.
\end{equation}
Consider $x\in F$ and the ball $B(x,r)$ with $r\in(0,1).$ Since
$\overline\Psi$ is a (LE) vector, there exists an $M>0$ such that
the number $N(x,r)$ of $\triangle^{(j)}$ (in the Moran cover of $F$)
with $\triangle^{(j)}\cap B(x,r)\neq\emptyset$ is bounded by $M$.
Hence
\begin{equation}
m_{\Psi}(B(x,r))\leq\sum_{j=1}^{N(x,r)}m_{\Psi}(\triangle^{(j)})
\leq\sum_{j=1}^{N(x,r)}L_{1}(\Psi^{(k)}_{\omega})^{\beta},
\end{equation}
where in the above summation $\omega$ corresponds to those for which
$\triangle_{\omega}=\triangle^{(j)}$, and $\triangle^{(j)}\cap
B(x,r)\neq\emptyset$. By
\eqref{equ_inf} and using $N(x,r)\leq M$, we have
\begin{equation}
m_{\Psi}(B(x,r))\leq
L_1M(\Psi^{(k)}_{\omega})^{\beta}\leq\frac{L_1M}{c_{*}^{\beta}}
(\Psi^{(k+1)}_{\omega})^{\beta}\leq\frac{L_1M}{c_{*}^{\beta}}r^{\beta}.
\end{equation}
This proves equation \eqref{eq.meas.reg}.

\hfill $\Box$

\begin{lemma}\label{lemma_hausdorff1}
Consider a MSC with a GMS, $F$. Suppose that (A3) holds,
$\overline\Psi=\{\overline\Psi^{(k)}\}$ is a \emph{basic} sequence
of vectors which satisfy the (UE), (LE) properties and
$c_{*}>0$. Then
$$\dim_{H}F=\mathrm{dim}_{H}(m_{\Psi})= s_{*}.$$
\end{lemma}

\noindent{\em Proof:} We first show that $\dim_{H}F\geq s_{*}$. Recall
$$\dim_{H}(m_{\Psi})=\inf\{\dim_H(E):
~\textrm{with}~m_{\Psi}(E)=1\}.$$ Moreover, from the proof of Lemma
\ref{lemma1}, we have $m_{\Psi}(B(x,r))\leq
\frac{L_1M}{c_{*}^{\beta}}r^{\beta}$, where $\beta<s_*$ is arbitrary. Hence
it follows that $s_{*}\leq\dim_{H}(m_{\Psi})\leq\dim_{H}F$, since
$\beta$ can be chosen arbitrarily close to $s_*$.

We now show that $\dim_H(F)\leq s_{*}$. Choose any $\beta>s_*$ then
for Hausdorff measure $\mathcal{H}^{\beta}$ we have
\begin{equation}\label{eq:upperh}
\begin{split}
\mathcal{H}^{\beta}(F)
&\leq\liminf_{k\rightarrow\infty}\sum_{\omega\in D_k}
\mathrm{diam}(\triangle_{\omega})^{\beta}\leq\liminf_{k\rightarrow\infty}\sum_{\omega\in
D_k}
C(\Psi^{(k)}_{\omega})^{\beta}\\
&\leq\liminf_{k\rightarrow\infty}\sum_{\omega\in D_k}
C(\Psi^{(k)}_{\omega})^{s_k}\leq
 CL_1\left(\liminf_{k\rightarrow\infty}\sum_{\omega\in D_k}
m_{\Psi}(\Delta_{\omega})\right)<\infty,
\end{split}
\end{equation}
where in the second line we take the infimum along the subsequence
$s_k$ such that $s_k<\beta$ (which holds infinitely often). Thus
$\mathcal{H}^{\beta}(F)\leq C$ and so $\dim_{H}F\leq\beta$. Since
$\beta>s_{*}$ is arbitrary, it follows that $\dim_{H}(F)\leq s_{*}$.
This completes the proof.

\hfill $\Box$


\begin{lemma}\label{lemma:box}
Consider a MSC with a GMS, $F$. Suppose that
$\overline\Psi=\{\overline\Psi^{(k)}\}$ is a basic collection of
vectors which satisfies the (UE), (LE) properties and $c_{*}>0.$
Then
$$\dim_{P}F\leq\overline{\dim_{B}}F\leq s^{*}.$$
\end{lemma}

\noindent{\em Proof:} We extend the ideas used in \cite[Page 141]{PesinWeiss}. Suppose (by contradiction) that
$\overline{\dim_{B}}(F)>s^{*}.$ Given the sequence $s_k$ and the
fact $s^{*}=\limsup_{k\to\infty}s_{k}$, then for all $\delta>0$,
there exists $\tilde{k}>0$ such that such that $\forall\,k\geq
\tilde{k}$, $\overline{\dim_{B}}(F)-3\delta>s_{k}$. By definition of
the pre-dimension sequence $s_{k}$ and noting that and corresponding
sequence of pressure functions $P_{k}:s\mapsto
P_{k}(s\log\Phi^{(k)})$ are decreasing in $s$, we have for all
$k\geq\tilde{k}$
\begin{equation}\label{equ_boxcontradiciton}
    P_{k}((\overline{\dim_{B}}(F)-3\delta)\log\Phi^{(k)})<0.
\end{equation}

Let $\beta=\overline{\dim_{B}}F$, then by the definition of upper
box dimension we have
$$\limsup_{\epsilon\rightarrow0}\frac{\log
N_{\epsilon}(F)}{-\log\epsilon}=\beta.$$
Hence given $\delta>0$, there is a sequence $\epsilon_n
=\epsilon_n(\delta)\to 0$, ($n\to\infty$), such that
\begin{equation}
N_{\epsilon}(F)\geq\epsilon_{n}^{\delta-\beta}.
\end{equation}

Given $\delta>0$, let $\epsilon$ be a representative from the sequence
$\epsilon_n$, which can be made arbitrarily small.
Let $\{\Delta^{j}\},j=1,\cdots,N^{\epsilon}(F)$ be the Moran covering of
$F$ at this $\epsilon$-scale. We have $N^{\epsilon}(F)\geq N_{\epsilon}(F).$
Since $0<d<1$ there exists $A>0$ such that for $j=1,\cdots,N^{\epsilon}(F):$
\begin{equation}\label{moran_cover1}
\frac{\epsilon}{A}\leq\Psi_{\omega}^{(n(\omega))} \leq\epsilon.
\end{equation}
Hence there exist uniform constants $C_{1}$ and $C_{2}$ such that
\[
C_{1}\log(\frac{1}{\epsilon})\leq n(x_{j})\leq
C_{2}\log(\frac{A}{\epsilon}).
\]
In the Moran covering the $n(\omega)$ can take on at most
$C_{3}:=C_{2}\log(\frac{A}{\epsilon})-C_{1}\log(\frac{1}{\epsilon})>0$
possible values. By the pigeon hole principle there
exists a positive integer $\alpha:=\alpha(\delta)$ with
$C_{1}\log(\frac{1}{\epsilon})\leq\alpha\leq
C_{2}\log(\frac{A}{\epsilon})$ such that for a sufficient small
$\epsilon$,
\begin{equation}\label{eq:box-estimate}
\sharp\{\omega~\mbox{such
that}~n(\omega)=\alpha\}\geq\frac{N^{\epsilon}(F)}{C_{3}}\geq\frac{N_{\epsilon}(F)}{C_{3}}
\geq\frac{\epsilon^{\delta-\beta}}{C_{3}}\geq\epsilon^{2\delta-\beta}.
\end{equation}
Recall that for any fixed number $s$, the potential $\Phi$, given by the function
$\Phi(x):=s\log\Psi_{\omega_{1}}(\mathbf{w}),$ where
$\mathbf{w}:=(\omega_1,\omega_2,\ldots)\in [D_{k}]$
is only dependent on the first coordinate $\omega_{1}$.
Therefore,
\[
(S_{n}\Phi)(x)=\sum_{j=1}^{n}\Phi(\sigma_{k}^{j}x)=s\log\prod_{j=1}^{n}\Psi^{(k)}_{\omega_j},
\]
and hence
$\exp(S_{n}\Phi)(x)=(\prod_{j=1}^{n}\Psi^{(k)}_{\omega_j})^{s}.$ We
have
\begin{equation}
P_{k}(\Phi(x))=\lim_{n\rightarrow\infty}\frac{1}{n}\log\sum_{(\omega_1,\ldots,\omega_n)}
\underset{x\in\Delta_{(\omega_1,\ldots,\omega_n)}}{\inf}\exp\left(\sum_{j=0}^{n-1}s\log([\sigma_{k}^{j}\mathbf{w}(x)]^{(1)})\right),
\end{equation}
and so
\begin{equation*}
P_{k}(\Phi(x))
=\log\left(\sum_{\omega\in D_k}(\Psi_{\omega}^{(k)})^{s}\right).
\end{equation*}
If we put $k=\alpha$ and apply equation \eqref{eq:box-estimate} with
$s=\overline{\dim_{B}}(F)-3\delta=\beta-3\delta$, we obtain
\begin{equation}\label{equ_countsequence}
\begin{split}
P_{\alpha}((\beta-3\delta)\Psi_{\omega_{1}}^{(\alpha)})
&\geq\log\left(\sum_{\omega\in\{\Delta^{(j)}\}}(\Psi_{\omega}^{(\alpha)})^{\beta-3\delta}\right)\\
&\geq \log\left(\left(\frac{\epsilon}{A}\right)^{\beta-3\delta}{\epsilon^{2\delta-\beta}}\right)\\
&=\log(\epsilon^{-\delta}A^{3\delta-\beta})\geq0.
\end{split}
\end{equation}
where $\{\Delta^{(j)}\}$ is the corresponding Moran cover. In the final
inequality we have used the fact that $A$ is independent of $\delta$.
Hence $\beta-3\delta<s_{\alpha}$. The constant $\alpha$ depends on the sequence
$\epsilon_n$, and can be taken arbitrarily large. This
implies that there is a subsequence $k_j\to\infty$, such that
$\beta-3\delta<s_{k_j}$, and hence $\beta-3\delta<\limsup s_k=s^*$ for every
$\delta>0$. The sequence $k_j$ implicitly depends on $\delta$, but for each
$\delta>0$ such an (infinite) sequence will always exist. Hence
$\beta\leq s^*$ in contradiction to \eqref{equ_boxcontradiciton}.

\hfill $\Box$

We now provide a lower bound for the box/packing dimensions when the
construction is conformal.
\begin{lemma}\label{lemma:conformalbox}
Suppose a GMS, $F$ arising from a MSC satisfies the conformal
condition as in Definition \ref{defconformal}, then
$$\overline{\dim_{B}}F= s^{*}.$$
\end{lemma}

\noindent{\em Proof:} Using Lemma \ref{lemma:box}, it suffices to
prove the lower bound. For each $\beta<s^{*}$, there exists a
subsequence $\{s_{k}\}$ such that for each $k$, $\beta<s_{k}$.
Moreover, the conformal condition implies that
$B(x,C^{-1}\Psi^{(k)}_{\omega})\subseteq\triangle_{\omega},$ for
each $\omega\in D_{k}.$  Using an equivalent definition of box dimension,
see Section \ref{sec:appendix}, we let
\[
W^{s}(F):=\lim_{r\to0}\sup\{\sum_{i}\operatorname{diam}(B_{i})^{s}:\operatorname{diam}(B_{i})\leq
r,B^{\circ}_{i}\cap B^{\circ}_{j}=\emptyset(i\neq j), B_{i}\cap
F\neq\emptyset\}.
\]
Then
\[
\overline{\operatorname{dim}_{B}}F:=\sup\{s:W^{s}(F)=\infty\}=\inf\{s:W^{s}(F)=0\},
\]
and hence,
\begin{equation*}
\begin{split}
W^{\beta}(F) &\geq\limsup_{k\to\infty}\sum_{\omega\in
D_{k}}\operatorname{diam}(B(x,C^{-1}\Psi^{(k)}_{\omega}))^{\beta}\\
&\geq\limsup_{k\to\infty}\sum_{\omega\in D_{k}} C^{-1}(\Psi^{(k)}_{\omega})^{s_{k}}\\
&\geq\limsup_{k\to\infty}\sum_{\omega\in
D_{k}}C^{-1}L_{1}^{-1}m_{\Psi}(\triangle_{\omega})>0,
\end{split}
\end{equation*}
which implies that $\overline{\dim_{B}}F>\beta$. Since $\beta$ is
arbitrary, it follows that
$\overline{\dim_{B}}=s^{*}.$

\hfill $\Box$

\begin{lemma}\label{lemma:conformalpacking}
Suppose a GMS, $F$ satisfies the conformal condition, then
$\dim_{P}F=\overline{dim_{B}}F.$
\end{lemma}
\noindent{\em Proof:} By Lemma \ref{lem_regular} in Section \ref{sec:appendix}, it suffices to show
that for any open set $V$, $\overline{\dim_{B}}(F\cap
V)\leq\overline{\dim_{B}}F=s^{*}$, provided $F\cap V\neq\emptyset$.
We do this as follows. Clearly $\overline{\dim_{B}}(F\cap
V)\leq\overline{\dim_{B}}F.$ Moreover, for any open set $V$ with
$F\cap V\neq\emptyset,$ there exists $\tilde\omega\in D_N$ such
that $\triangle_{\tilde{\omega}}\subset V$. Taking $\beta<s^{*}$, there
exists a subsequence $s_{k}$ with $\beta<s_{k},$ such that
\begin{equation*}
\begin{split}
W^{\beta}(F\cap V) &\geq\limsup_{k\to\infty}\sum_{\omega\in
D_{k},\triangle_{\omega}\subseteq\triangle_{\tilde\omega}}\operatorname{diam}(B(x,C^{-1}\Psi^{(k)}))^{\beta}\\
&\geq\limsup_{k\to\infty}\sum_{\omega\in
D_{k},\triangle_{\omega}\subseteq\triangle_{\tilde\omega}}
C^{-1}(\Psi^{(k)})^{s_{k}}\\
&\geq \limsup_{k\to\infty}\sum_{\omega\in
D_{k},\triangle_{\omega}\subseteq\triangle_{\tilde\omega}}
C^{-1}L_{1}^{-1}m_{\Psi}(\triangle_{\omega})\\
&=C^{-1}L_{1}^{-1}m_{\Psi}(\triangle_{\tilde{\omega}})>0.
\end{split}
\end{equation*}
Hence we obtain $s^{*}=\overline{\dim_{B}}(F\cap
V)=\overline{\dim_{B}}(F),$ which completes the proof.

\hfill $\Box$

\subsection{Proof of Theorem \ref{main.theorem.2}}

The proof of Theorem \ref{main.theorem.2} is as follows. We claim
first of all that $\dim_{H}F\leq s_{*}$. The proof of this claim
follows step by step the proof of Lemma \ref{lemma_hausdorff1} via
equation \eqref{eq:upperh}. Hence it suffices to show only that
$\dim_{H}(F)\geq s_{*}$.

Suppose $\beta<s_*$, then there exists a $K\in\mathbb{N}$ such that
for all $k\geq K$, $s_{k}>\beta$. Moreover for any $\omega\in D$,
$(\Psi^{(k)}_{\omega})^{s_k}<(\Psi^{(k)}_{\omega})^{\beta}.$

Since $\log d_k/\log M_k\to 0$, there exists an $\epsilon>0$ such
that for all $k\geq K$, $M^{\epsilon/2}_k<d^{\beta+\epsilon}_k$ and
hence that $M^{\epsilon/2}_k/d^{\beta+\epsilon}_k<1$.

Now take the Moran cover $\triangle^{(j)}$ such that
$m_{\Psi}(\triangle^{(j)})\leq L_1(\Psi^{(k)}_{\omega})^{s_k}
<(\Psi^{(k)}_{\omega})^{\beta}.$ Given $r>0$, the Moran cover
$\triangle^{(j)}$ has the property that
$(\Psi^{(n(\omega))}_{\omega})\geq r$ and
$(\Psi^{(n(\omega)+1)}_{\omega})< r$. Since $s_{k}>\beta$ we choose
$\epsilon$ sufficiently small so that $s_{k}>\beta+\epsilon.$
Therefore, we obtain the following series of estimates:
\begin{equation}
\begin{split}
m_{\Psi}(B(x,r)) &\leq\sum_{j=1}^{N(x,r)}m_{\Psi}(\triangle^{(j)})
\leq\sum_{j=1}^{N(x,r)}L_{1}(\Psi^{(k(\omega))}_{\omega})^{\beta+\epsilon}\\
&\leq
\sum_{j=1}^{N(x,r)}\frac{L_1}{d^{\beta+\epsilon}_{k(\omega)}}(\Psi^{(k(\omega)+1)}_{\omega})^{\beta+\epsilon}
\leq \sum_{j=1}^{N(x,r)}\frac{L_1M^{\epsilon/2}}{d^{\beta+\epsilon}_{k(\omega)}}(\Psi^{(k(\omega)+1)}_{\omega})^{\beta+\frac{\epsilon}{2}}\\
&\leq L_1Mr^{\beta+\frac{\epsilon}{2}}.
\end{split}
\end{equation}
It follows that
$\dim_{H}(m_{\Psi})\geq s_{*}$, since we can choose $\beta$
arbitrarily close to $s_*$ and $\epsilon$ arbitrarily close to 0. It
follows that $\dim_{H}(F)\geq \dim_H(m_{\Psi})\geq s_{*}$ and hence
we obtain $\dim_{H}(F)=s_{*}$

\medskip

We now turn to the box dimension. First consider the case where the
vector $\overline\Psi$ is upper-estimating, but the construction is
not conformal. We can repeat the proof of Lemma \ref{lemma:box}, but
we note that the constant $A$ appearing in equation
\eqref{moran_cover1} is now dependent on $\epsilon$.  Taking again
the Moran covering
$\{\Delta^{j}_{\omega}\},j=1,\cdots,N^{\epsilon}(F)$ of $F$ at scale
$\epsilon$, we have $\Psi^{(n(\omega)+1)}_{\omega}<r\leq
\Psi^{(n(\omega))}_{\omega}$. Recalling that $d_k=\min\limits_{1\leq
j\leq n_{k+1}}c_{j}^{(k+1)}$, we obtain
$$\epsilon>\Psi^{(n(\omega)+1)}_{\omega}=\Psi^{(n(\omega))}_{\omega}c^{(n+1)}_j
(\omega)\geq d_{n(\omega)+1}\Psi^{(n(\omega))}_{\omega},$$ and hence
$\Psi^{(n(\omega))}_{\omega}\leq\epsilon d^{-1}_{n(\omega)+1}.$
Since $\lim_{k\to\infty}\frac{\log d_{k}}{\log M_{k}}=0$ it follows
that for all $\eta>0$, there exists a $K$, such that $\forall\,k\geq
K$, $1>d_k>M_{k}^{\eta}>0$, and therefore
$$\Psi^{(n(\omega))}_{\omega}\leq\epsilon M^{-\eta}_{k}.$$
Hence by definition of $M_k$ we obtain (for arbitrary $\eta>0$):
$\epsilon<\Psi^{n(\omega)}_{\omega}\leq\epsilon^{\frac{1}{1+\eta}}.$
Following step by step the proof of Lemma \ref{lemma:box}, we obtain
$\overline{\dim}_{B}(F)\leq s^*.$

\medskip

When the construction is conformal, the upper bound for $\dim_{P}F$ and $\overline{\dim_{B}}F$
is obtained as in the calculation directly above. The lower bounds follow from Lemmas~\ref{lemma:conformalbox} and
\ref{lemma:conformalpacking}.

\hfill $\Box$

\subsection{Proof of Theorem \ref{main.theorem3}}
To prove Theorem \ref{main.theorem3} we consider a truncated construction, such
as that considered in \cite{Hua00}. We
remove words $\omega\in D_k$ for which
$c^{(k)}_j<\epsilon$, and define
\begin{equation*}
D_{k}(\epsilon)=\{\omega\in
D_k: c^{(k)}_j\geq\epsilon,\forall i\leq
k-1\},\quad \tilde{D}_k(\epsilon)=\{\triangle_{\omega}:\omega\in
D_k(\epsilon)\},
\end{equation*}
and
\begin{equation*}
E_{k}(\epsilon)=\bigcup_{\omega\in
D_{k}(\epsilon)}\triangle_{\omega},\qquad F(\epsilon)=\bigcap_{k\geq
0} E_{k}(\epsilon).
\end{equation*}
For the $\epsilon$-truncated construction $F(\epsilon)$ of $F$, the
associated vectors $\{\Psi^{(k)}_{\omega}\}$ are both upper and
lower-estimating, and so we can use Theorem \ref{main.theorem} to
find the fractal dimension of $F(\epsilon).$ In particular the
dimension of $F(\epsilon)$ can be found by taking appropriate limits
along the pre-dimension sequences $s_{k}(\epsilon)$, where
$s=s_{k}(\epsilon)$ solves the equation
$P_{k}(s\log(\Psi^{(k)}_{\omega}\mathcal{X}_{D_{k}(\epsilon)}(\omega)))=0$.
Here $\mathcal{X}_{D_{k}(\epsilon)}(\omega)$ denotes the indicator
function of $D_{k}(\epsilon)$. The following lemma makes explicit
the relation between $s_{k}(\epsilon)$ and $s_{k}$, the latter value
being the solution to $P_{k}(s\log(\Psi^{(k)}))=0$.
\begin{lemma}\label{lem:inf-zero}
Suppose $s_k$ and $s_k(\epsilon)$ are solutions to the respective
pressure equations
$$P_{k}(s\log(\Psi^{(k)}_{\omega}))=0,\quad
P_{k}(s\log(\Psi^{(k)}_{\omega}\mathcal{X}_{D_{k}(\epsilon)}(\omega)))=0.$$
Suppose that $s_*=\liminf s_k>0$. Then for $k$ sufficiently large,
\begin{equation}
0\leq s_k-s_k(\epsilon)\leq\mathcal{O}(\epsilon^{s*/2}),
\end{equation}
where the implied constant in $\mathcal{O}(\cdot)$ is independent of
$k$.
\end{lemma}

\noindent{\em Proof:}
The arguments follow close to \cite{Hua00} and we provide the main steps.
Suppose that $s_*>0$. Firstly, there exist constants $0<\alpha,\beta<1$ such that
\begin{equation*}
 \alpha<\inf_{k}\max_{1\leq j\leq n_k}c^{(k)}_j,\quad
\sup_{k}\max_{1\leq j\leq
n_k}c^{(k)}_j <\beta.
\end{equation*}
Observe that for any $\omega\in D_k$ (and $i\leq k$) we have that
\begin{equation}\label{infzero-est1}
\sum_{j=1}^{n_i}(c^{(i)}_j)^{s_k}\geq\alpha^{s_k}\geq\alpha^{\tilde{d}},
\end{equation}
\begin{equation}\label{infzero-est2}
\sum_{j=1}^{n_i}(c^{(i)}_j)^{s_k}\mathcal{X}_{D_{k}(\epsilon)}(\omega)>
\sum_{j=1}^{n_i}(c^{(i)}_j)^{s_k}-M\epsilon^{s_k}.
\end{equation}
where $\tilde{d}$ is the dimension of the space. From equation \eqref{infzero-est2} we obtain
\begin{equation}\label{infzero-est3}
\sum_{j=1}^{n_i}(c^{(i)}_j)^{s_k}\mathcal{X}_{D_{k}(\epsilon)}(\omega)>
\left(1-\frac{M\epsilon^{s_k}}{\alpha^{\tilde{d}}}\right)\sum_{j=1}^{n_i}(c^{(i)}_j)^{s_k}.
\end{equation}
Now, for any $\gamma>0$ we have
\begin{equation}\label{infzero-est4}
\sum_{j=1}^{n_i}(c^{(i)}_j)^{s_k}\mathcal{X}_{D_{k}(\epsilon)}(\omega)=
\sum_{j=1}^{n_i}(c^{(i)}_j)^{s_k-\gamma}(c^{(k)}_j)^{\gamma}\mathcal{X}_{D_{k}(\epsilon)}(\omega)
\leq\beta^{\gamma}\sum_{j=1}^{n_i}(c^{(k)}_i)^{s_k-\gamma}\mathcal{X}_{D_{k}(\epsilon)}(\omega).
\end{equation}
Taking products and combining equations \eqref{infzero-est3}, \eqref{infzero-est4} we obtain
\begin{equation}\label{infzero-est5}
\prod_{i=1}^{k}\sum_{j=1}^{n_i}(c^{(i)}_j)^{s_k-\gamma}
\mathcal{X}_{D_{k}(\epsilon)}(\omega)\geq
\beta^{-k\gamma}\left(1-\frac{M\epsilon^{s_k}}{\alpha^{\tilde{d}}}\right)^k
\prod_{i=1}^{k}\sum_{j=1}^{n_i}(c^{(i)}_j)^{s_k}.
\end{equation}
We observe that $s_k$ and $s_{k}(\epsilon)$ satisfy the pre-dimension equations
\begin{equation}\label{infzero-est6}
\prod_{i=1}^{k}\sum_{j=1}^{n_i}(c^{(i)}_j)^{s_k(\epsilon)}
\mathcal{X}_{D_{k}(\epsilon)}(\omega)=1,\quad\prod_{i=1}^{k}\sum_{j=1}^{n_i}(c^{(i)}_j)^{s_k}=1.
\end{equation}
Since $s_*>0$, there exists $k_0$, such that $s_k>s_*/2$ for all $k\geq k_0$. Moreover, there exists $\epsilon_0$,
such that for all $\epsilon<\epsilon_0$, we have
\begin{equation}
0<\gamma_{\epsilon}:=\log(1-M\alpha^{-\tilde{d}}\epsilon^{s_k})/(\log\beta)
<\log(1-M\alpha^{-\tilde{d}}\epsilon^{\frac{s_*}{2}})/(\log\beta)<\frac{s_*}{2}.
\end{equation}
From equations \eqref{infzero-est5} and \eqref{infzero-est5}, we see that for any $\gamma<\gamma_{\epsilon}$ we
have
\begin{equation}
\prod_{i=1}^{k}\sum_{j=1}^{n_i}(c^{(i)}_j)^{s_k-\gamma}
\mathcal{X}_{D_{k}(\epsilon)}(\omega)\geq 1,
\end{equation}
and therefore have $s_k(\epsilon)>s_k-\gamma$. From the observation that $s_k\geq s_{k}(\epsilon)$ the result now follows.

\hfill $\Box$

We now claim that $\dim_H(F)=\liminf s_{k}=s_*$. Since
$F(\epsilon)\subset F$, we have $\dim_H(F)\geq s_*(\epsilon)$, and
by Lemma \ref{lem:inf-zero}, we have $\lim_{\epsilon\to
0}s_*(\epsilon)=s_*$. Hence $\dim_H(F)\geq s_*$. For the upper bound
we just apply the same argument as Lemma \ref{lemma_hausdorff1}.

For the upper-box and packing dimensions we claim that
$\overline{\dim}_B(F)=\limsup s_k=s^*$.
For a monotonically
decreasing sequence $\epsilon_n\to 0$ let
$F^{*}=\bigcup_{n=1}^{\infty}F(\epsilon_n)$. Then by the closure property of upper-box
dimension we have $\overline{F^*}=F$, and so
$\overline{\dim}_B(F^*)=\overline{\dim}_B(F)$. It therefore suffices to
calculate $\overline{\dim}_B(F^*)$. By Theorem \ref{main.theorem}, and
$\forall~\epsilon>0$, we have
$\overline{\dim}_BF(\epsilon)=\dim_{P}(F(\epsilon))$. These dimensions equal
$\limsup_k s_k(\epsilon).$ Furthermore, $\dim_{P}(F^*)=\limsup_{n\to\infty}\dim_P
(F(\epsilon_n))=s^{*}.$ This completes the proof.

\section{Appendix: Background on fractal dimension and its computation}\label{sec:appendix}
\subsection{Hausdorff, Box and Packing dimensions}

In this section we give the relevant background on dimension theory, see \cite{Falconer97, Falconer, Mattila95}
for a more general discussion.

Suppose $F$ is a non-empty subset in $\mathbb{R}^{d}$. For any non-negative number $s$ and
$\epsilon>0$, let
\begin{equation}\label{equ_hausdorff}
\mathcal{H}^{s}_{\epsilon}(F):=
\inf\left(\sum_{i}(\mathrm{diam}(U_i))^{s}\right),
\end{equation}
where the infimum is taken over all covers $\{U_{i}\}$ with
$\mathrm{diam}(U_i)<\epsilon$. As $\epsilon$ decreases, the class of
permissible covers of $F$ in \eqref{equ_hausdorff} is reduced, and
therefore, the infimum $\mathcal{H}^{s}_{\epsilon}$ increases. The
limit $\mathcal{H}^{s}(F):=\lim_{\epsilon\to
0}\mathcal{H}^{s}_{\epsilon}(F)$ exists, and is called as
\emph{Hausdorff measure}. The corresponding \emph{Hausdorff dimension} of $F$ is
defined by
\begin{equation}\label{equ_hausdorff dimension}
\mathrm{dim}_{H}(F):=\inf\{s\geq0:\mathcal{H}^{s}(F)=0\}=\sup\{s\geq0:\mathcal{H}^{s}(F)=\infty\}.
\end{equation}
A disadvantage of Hausdorff dimension lies in its calculation. Alternative definitions of dimension, which are
perhaps easier to estimate are the following.

The (upper) box dimension is relatively easier to estimate than Hausdorff dimension, and is defined as follows. Given a
non-empty set $F\subset\mathbb{R}^{d}$, and $\epsilon>0$, let
$N(\epsilon)$ denote the smallest number of $\epsilon$-balls needed
to cover $F$. The (upper) box dimension of $F$ is defined by:
\begin{equation}\label{equ_boxdimension}
\overline{\mathrm{dim}}_{B}(F):=\lim\sup_{\epsilon\to 0}\frac{\log
N(\epsilon)}{\log(1/\epsilon)}.
\end{equation}
Analogous to Hausdorff dimension, there is an alternative
description of upper box dimension \cite{Pe97}: for any non-negative
number $s$, let
\begin{equation}\label{equ_quasimeasure}
\begin{array}{rl}
W^{s}(F):=\lim_{\epsilon\to0}\sup\{\sum_{i}\operatorname{diam}(B_{i})^{s}:&
B_{i}~\mbox{is a ball with}~\operatorname{diam}(B_{i})\leq
\epsilon,\\
&B^{\circ}_{i}\cap B^{\circ}_{j}=\emptyset(i\neq j), B_{i}\cap
F\neq\emptyset\},
\end{array}
\end{equation}
then
\begin{equation}\label{equ_boxanalogs}
\overline{\operatorname{dim}_{B}}(F):=\sup\left\{s:W^{s}(F)=\infty\}=\inf\{s:W^{s}(F)=0\right\}.
\end{equation}
It is worth mentioning that $W^{s}(\cdot)$ in equation
\eqref{equ_quasimeasure} usually does not define a measure (due to
lack of subadditivity). Moreover $\overline{\operatorname{dim}_{B}}(F)=\overline{\operatorname{dim}_{B}}(\overline {F})$,
where $\overline{F}$ is the closure of $F$. Hence box dimension can give positive values to countable sets (unlike Hausdorff dimension).

Comparing equation \eqref{equ_hausdorff} with equation
\eqref{equ_quasimeasure}, we have
\begin{equation}\label{equ_gapofdiemensions}
\mathrm{dim}_{H}(F)\leq \overline{\mathrm{dim}}_{B}(F).
\end{equation}

We now introduce packing dimension and packing measure. Let
$$\mathcal{P}^{s}_{\epsilon}(F):=\sup\{\sum_{i}\operatorname{diam}(B_i)^s\}$$
where the supremum is taken over a collection of disjoint balls
$\{B_{i}\}$ of radius at most $\epsilon$ and with centers in $F$.
The limit $\mathcal{P}^{s}_{0}(F):=\lim_{\epsilon\to
0}\mathcal{P}^{s}_{\epsilon}(F)$ exists. However, by considering
countable dense sets, it is easy to see that $\mathcal{P}^{s}_{0}$
is not a measure (again, due to lack of subadditivity). Hence, we
modify $\mathcal{P}^{s}_{0}$ to
\begin{equation}\label{equ_packing dimension}
\mathcal{P}^{s}(F):=\inf\left\{\sum_i\mathcal{P}^{s}_{0}(F_i):
F\subset \bigcup_{i=1}^{\infty}F_i\right\},
\end{equation}
which is a measure, and is called an $s$-dimensional packing
measure. The \emph{packing dimension} is naturally defined as
\begin{equation}
\mathrm{dim}_{P}(F):=\sup\{s:\mathcal{P}^{s}(F)=\infty\}=\inf\{s:\mathcal{P}^{s}(F)=0\}.
\end{equation}

For a general set $F\subset\mathbb{R}^{d}$, the following relations
hold:
\begin{equation}\label{equ_fractdimensionorder}
\mathrm{dim}_{H}(F)\leq \mathrm{dim}_{P}(F)\leq
\overline{\mathrm{dim}}_{B}(F), ~~\mbox{and}~~
\mathcal{H}^{s}(F)\leq\mathcal{P}^{s}(F).
\end{equation}
Suitable examples show that none of inequalities in
\eqref{equ_fractdimensionorder} can be replaced by equalities
\cite{Falconer}.

The following lemma is useful for studying packing and box
dimension, especially for fractal sets with some degree of self
similarity.
\begin{lemma}\cite[Corollary 3.9]{Falconer97}\label{lem_regular}
Let $F\subset\mathbb{R}^{n}$ be compact and for all open sets $V$ that intersect with $F$ suppose that
$\overline{\dim_{B}}(F\cap V)=\overline{\dim_{B}}(F).$
Then
$\dim_{P}(F)=\overline{\dim_{B}}(F).$
\end{lemma}

\subsection{Relation between conformality and the (LE) property}
The following lemma gives the relationship between a conformal construction and a construction which admits
a lower estimating vector.
\begin{lemma}\label{Conformality implies l-estimating} If a vector $\overline\Psi$ is conformal, and
$c_{*}>0$,
then the vector $\overline\Psi$ satisfies (LE) property.
\end{lemma}

\noindent{\em Proof:} For any fixed $0<r<1,$ and any $x\in F$, consider
the open ball $B(x,r)$ centered in $x$ with radius of $r$, and let
$N(x,r)$ be the number of Moran covering $\{\triangle^{(j)}\}$ that
have nonempty intersection with $B(x,r).$ Hence
\[
B(x,r)\bigcup\left(\cup_{j=1}^{N(x,r)}\triangle^{(j)}\right)\subseteq
B(x,R),
\]
where
\[
R=2r+\sup_{j}\operatorname{diam}(\triangle^{(j)}).
\]
Using the conformal condition and recalling from the definition of
$\triangle^{(j)}$, we can choose elements $\triangle_{\omega}$,
$\omega\in D_{n+1}$
such that $\Psi^{(n)}_{\omega}\geq r$ and $\Psi^{(n+1)}_{\omega}\leq
r$. Therefore,
\[
R\leq 2r+\sup_{j}C\Psi^{(k)}_{\omega_j}\leq 2r+\frac{r}{c_*},
\]
and
\[
\operatorname{diam}(\triangle^{(j)})\geq
C^{-1}\Psi^{(k)}_{\omega_j}\geq C^{-1}r.
\]
Hence it follows that for each $x\in F$ and $0<r<1$
\[
N(x,r)\leq\frac{2r+c_{*}^{-1}r}{C^{-1}r}=\frac{2+c_{*}^{-1}}{C^{-1}}<\infty.
\]
Therefore the vector $\overline\Psi$ satisfies (LE) property.

\hfill $\Box$

\subsection{Background on thermodynamic formalism}
For inhomogeneous Moran set constructions we used intermediate constructions based on finite symbolic schemes to calculate the
fractal dimension. We review relevant background on thermodynamic formalism for these
finite symbolic schemes, see for example \cite{Falconer97, PP90, PesinWeiss}.

Consider the finite symbolic dynamical system $(\sum_{p}^{+},\sigma)$, where
$\sum_{p}^{+}=\{0,\cdots,p-1\}^{\mathbb{N}}$ and $\sigma:\sum_{p}^{+}\to\sum_{p}^{+}$ as the left-shift map.
Suppose $Q\subset\Sigma^{+}_p$ is a $\sigma$-invariant set.
If $\omega\in Q$, then we write $\omega=(i_1,i_2,\ldots)$, with
$i_j\in\{0,\ldots, p\}$ an admissible sequence (for $j\geq 1$). We turn $\Sigma^{+}$ into a metric space using a
standard symbolic metric, such as that given in Section \ref{sec:definitions}.
Given $\omega\in Q$, we write $C_{i_1,\ldots i_k}(\omega)\subset Q$ as the $k$-length
cylinder set that contains $\omega$. Given an $\alpha-$H\"{o}lder continuous function
$\phi:Q\to\mathbb{R}^{+}$, let
$S_k(\phi):=\sum_{i=0}^{k-1}\phi\circ\sigma^{i}$, then the
\emph{(topological) pressure} $P(\phi)$ is defined by
\begin{equation}\label{def_press_fun}
P(\phi):=\lim_{n\rightarrow\infty}\frac{1}{n}\log\left(\sum_{\underset{\textrm{admissible}}{(i_1,\ldots, i_n)}}
\inf_{\omega\in C_{i_1,\ldots, i_n}}\exp(S_{n}(\phi)(\omega))\right).
\end{equation}
For topological dynamical systems, the following variational principle holds. Let $\mathcal{M}(Q)$ denote the space of
$\sigma$-invariant measures on $Q$. Then for $\phi:Q\to\mathbb{R}^{+}$ H\"older continuous we have
$$P(\phi)=\sup_{\mu\in\mathcal{M(Q)}}\left(h_{\mu}(\sigma)+\int_{Q}\phi\,d\mu\right),$$
where $h_{\mu}(\sigma)$ is the topological entropy of $\sigma$.
The measure $\mu=\mu_{\phi}$ that gives rise to the supremum is called an \emph{equilibrium measure}. This measure always exists,
but need not be unique. Another measure of significance is that of a \emph{Gibbs measure}. For any $\alpha-$H\"{o}lder continuous map
$\phi:\sum_{p}^{+}\to\mathbb{R}^+$, an invariant measure $\mu$ is called a Gibbs measure for the potential
$\phi$ if there exists a constant $D>1$ such that
\begin{equation}\label{equ_Gibbs measure}
   D^{-1}\leq\frac{\mu\{y:y_{i}=x_{i},
   i=1,\cdots,n\}}{\exp(-nP(\phi)+\sum_{k=0}^{n-1}\phi(\sigma^{k}(x)))}<D
\end{equation}
for all $x=(x_{1},x_{2},\cdots)\in\sum_{p}^{+}$ and $n\geq0.$ In
fact, for the shift map $\sigma$ on a finite symbolic space, the
hypothesis of the $\alpha-$H\"{o}lder continuity of the potential
$\phi$ ensures the existence and uniqueness of the Gibbs measure and
its coincidence with the equilibrium state for $\phi.$ However for
more general symbolic schemes less is known about the existence of
such measures. To study results on fractal dimension, the potential
$\phi$ of interest is that which depends only on the first
coordinate, i.e., $\phi(x)=\phi(x_{1})$. In \cite{PP90} it shown
that for given numbers $0<\lambda_{i}<1, i=1,\cdots,p,$ and
potential function $\phi:\sum_{p}^{+}\to\sum_{p}^{+}$ defined by
$\phi(x)=\phi(x_{1},x_{2},\cdots)=\log\lambda^{-1}_{x_{1}}$, the
equation $P(s\phi)=0$ has a unique solution in $s$. Moreover $\phi$
is H\"{o}lder continuous. This unique solution $s$ is equal to the
Hausdorff, Packing and Boxing dimensions of certain repelling
invariant sets generated by IFS, see for example \cite{Falconer97,
Pe97}.

\subsection*{Acknowledgement}
In preparing this work we would like to thank T. Jordan and J. Rivera-Letelier
for useful discussions.

%
%
%
%
%
%
%
%
%
%
%
%
%

\end{document}